\documentclass[a4paper,12 pt,oneside,reqno]{amsart}
\usepackage{graphicx}
\DeclareGraphicsExtensions{.pdf,.eps,.png,.jpg}
\usepackage{amssymb,amsthm,amsmath,mathrsfs,amscd} 
\usepackage{tikz,pgfplots}
\usetikzlibrary{intersections,calc}
\usepackage{geometry}
\usepackage{booktabs}
\usepackage{xcolor}
\usepackage{float}
\usepackage[pagewise]{lineno}

\newcommand{\eps}{\varepsilon}

\newcommand{\R} {\mathbb{R}}

\newtheorem{teorema}{Theorem}[section]

\newtheorem{proposizione}[teorema]{Proposition}

\theoremstyle{definition}

\theoremstyle{remark}

\usepackage[autostyle,italian=guillemets]{csquotes}
\usepackage[style=numeric,backend=biber,maxnames=10]{biblatex}
\usepackage{amsaddr}

\begin{document}
\title{Stability analysis of microscopic models for traffic flow with lane changing}
\author{Matteo Piu\textsuperscript{*}}
\thanks{\textsuperscript{*} Corresponding author.}
\address{Dipartimento di Scienze di Base e Applicate per l'Ingegneria,  Sapienza - Universit\`{a} di Roma}
\email{matteo.piu@uniroma1.it}
\author{Gabriella Puppo}
\address{Dipartimento di Matematica,  Sapienza - Universit\`{a} di Roma}
\email{gabriella.puppo@uniroma1.it}
\date{\today}

\begin{abstract}
 {This paper investigates the mathematical modeling and the stability of multi-lane traffic in the microscopic scale, studying a model based on two interaction terms.
To do this we propose  simple lane changing conditions and we study the stability of the steady states starting from the model in the one-lane case and extending the results to the generic multi-lane case with the careful design of the lane changing rules. We compare the results with numerical tests, that confirm the predictions of the linear stability analysis and also show that the model is able to reproduce stop \& go waves, a typical feature of congested traffic.}
\end{abstract}
\maketitle

\subsection*{AMS subject classifications.} 76A30, 34D20, 70-10.

\subsection*{Keywords.} Traffic flow, multi-lane traffic, lane changing, stability analysis.

\section{Introduction}
 {In this paper we deal with microscopic modeling of traffic flow, focusing on lane changing dynamics. In particular we study a second order model for one lane that combines two different interaction terms and we describe the extension to the  {multi-lane} case giving particular attention at the two-lane case.}

\subsection{Related work.} The interest in the dynamics of traffic flow dates back to the first half of the twentieth century and the related mathematical literature is quite large. An overall view can be found, for instance, in the book by Haberman \cite{MR1600250} and in the survey paper by Helbing \cite{Helbing20011067}.

There are various points of view for modeling traffic flow. In this paper we concentrate on the microscopic approach that is based on the dynamics of individual vehicles considering the individual behaviour of each driver. A typical microscopic model is the Car Following model or Follow the Leader model (FtL) based on the idea that the dynamics of each vehicle (follower) depends on the vehicle in front (leader) and therefore the other vehicles do not affect it. These models are normally for  {single-lane} roads \cite{Brackstone1999181,MR94251,MR0135615}. A typical Follow the Leader model can be described as follows. In a  {single-lane} with $N$ vehicles where overtaking is not allowed, we are interested in study the position $x_n(t)$ and the velocity $v_n(t)$ of each vehicle $n=1,\dots,N$ at different times $t$. This dynamics can be described by a system of ordinary differential equations:
\begin{equation}
\begin{cases}
\dot{x}_n(t)=v_n(t) &n=1,\dots,N\\
\dot{v}_n(t)=a(x_n(t),x_{n+1}(t),v_n(t),v_{n+1}(t))  &n=1,\dots,N-1\\
\dot{v}_N=w(t).
\end{cases}
\end{equation}
where $a(\cdot)$ is a given acceleration function and $w(\cdot)$ is the dynamics of the leader vehicle, independent from the other vehicles (followers).

Many  {single-lane} car following models have been developed and applied to study traffic dynamics. Here we recall some models that will be useful in the following.

The Follow the Leader model, introduced in \cite{Reuschel1950193,Reuschel2}, assumes that each vehicle modifies its velocity based on the distance (headway) $x_{n+1}-x_n$ to the vehicle ahead, the $n+1$-th, and to the difference in velocities between its own velocity $v_n$ and the velocity of the vehicle ahead $v_{n+1}$, multiplied by appropriate coefficients $\beta_n$. This model can be described by the following system
\begin{equation}\label{ftl_classico}
\begin{cases}
\dot{x}_n(t)=v_n\\
\dot{v}_n(t)= \beta_n \frac{v_{n+1}-v_n}{(x_{n+1}-x_n)^2} \\
\end{cases}.
\end{equation}

The optimal velocity model (OVM) of Bando et al. \cite{Bando19951035,Bando1994203} in which a driver aims to a desired velocity function $V$ that depends on the headway with the vehicle ahead. The equation of this model is given by
\begin{equation}\label{bando_classico}
\begin{cases}
\dot{x}_n(t)=v_n\\
\dot{v}_n(t)= \alpha_n (V(x_{n+1}-x_n)-v_n) \\
\end{cases}
\end{equation}
with appropriate coefficients $\alpha_n$. 

We mention also some interesting works. Pipes proposed \cite{Pipes1953274} a traffic model in which each vehicle maintains a certain prescribed ''following distance'' from the preceding vehicle; the generalized force model (GFM) by Helbing and Tilch \cite{Helbing1998133} in which the optimal velocity function is obtained calibrating the parameters with the observed data; the full velocity difference model (FVDM) by Jiang et al. \cite{Jiang2001017101/1} that predicts  delay time of car motion and kinematic wave speed at jam density; the optimal velocity difference model (OVDM) by Peng et al. \cite{Peng20113973} where  a new term is introduced involving the optimal velocity functions and the vehicles $n,n+1,n+2$. Aw et al. \cite{Aw2002259} studied the derivation of a continuum model starting from the FtL model.
We mention an analytical study for the OVM with a stepwise specification of the optimal velocity function and a simple kind of perturbation in \cite{helbingnuovo}. 

Another type of microscopic model is given by lane changing models which provide for the possibility of changing lanes according to the analysis of some factors that intervene in the decision process, for example the need, opportunity and safety of a lane change \cite{Errampalli,Toledo200765}. The interest in modeling vehicle lane changing is due to the effect that it induces in traffic flow, for instance in bottleneck discharge rate and in the  {stop \& go} oscillations. Here we recall some works. Cassidy and Rudjanakanoknad \cite{Cassidy2005896} showed that when traffic density upstream of a busy merge increases beyond a critical value, vehicles manoeuvre toward faster lanes causing traffic breakdown and ''capacity drop'' of the road; Zheng et al. \cite{Zheng20111378} showed that lane changing are responsible for transforming subtle localized oscillations to substantial disturbances; Klar and Wegener \cite{Klar19981002,Klar1998983} developed a model based on reaction thresholds from which they derived a kinetic model; Song and Karni \cite{Song20191429} proposed a macroscopic model in which the acceleration terms take lead from microscopic car-following models, and yield a non-linear hyperbolic system with viscous and relaxation terms; Herty et al. \cite{Herty20182252}  proposed a macroscopic model, which accounts for lane-changing on motorway, based on a two-dimensional extension of the Aw and Rascle and Zhang macroscopic model for traffic flow; Gong et al. \cite{Gong20211485} presented a  finite dimensional hybrid system  based on the continuous Bando-Follow-the-Leader dynamics coupled with discrete events due to lane changing; Goatin and Rossi \cite{Goatin20191967} developed a macroscopic model for  {multi-lane} road networks with discontinuities both in the speed law and in the number of lanes; Hodas and Jagota presented in \cite{Hodas20031247} a microscopic model for  {multi-lane} dynamics where each car experiences a force resulting from a combination of the desire of the driver to attain a certain velocity and change of the force due to cars interactions; Kesting et al. \cite{Kesting200786} proposed a general model to derive lane changing rules for discretionary and mandatory lane changes for a wide class of car following models; Lv et al. in \cite{Lv20112303} extended the continuous  {single-lane} models to simulate the lane changing behaviour on an urban roadway with three lanes  and in \cite{Lv20131142} proposed a model where lane changing is not instantaneous but is a continuing process which can affect the following cars; Zheng et al. in \cite{Zheng2013367} analysed the effects of lane changing in the driver behaviour.

\subsection{Goal and paper organization} This paper proposes the study of a second order microscopic model combining models \eqref{ftl_classico} and \eqref{bando_classico} for reproducing traffic flow and its extension to the  {multi-lane} case with simple lane changing conditions in order to study its stability under perturbations . In Section 2 we introduce the model for a  {single-lane} and we study its stability in the linearized case, then we show numerical tests making comparisons with model \eqref{bando_classico}. In section 3 we describe the extension of the model to the two-lane case studying its stability around the equilibrium when a lane is perturbed. We present some numerical tests that confirm the predictions of the linear stability analysis. Finally, in section 4, we illustrate the generalization of the model to the generic  {multi-lane} case.

\section{Single-lane model} 
\subsection{Description}
In this section we describe the main mathematical model we use in this paper. Consider a homogeneous population of $N \in \mathbb{N}$  vehicles, and denote by $x_n=x_n(t)$ and $v_n=v_n(t)$ the position and the velocity of the $n$-th vehicle at time $t\in \R^+$. We want to describe the traffic flow in a road with a  {single-lane} where overtaking is not allowed.

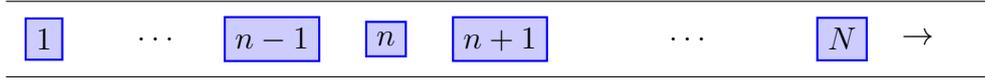
\begin{figure}[H]
\centering
\begin{tikzpicture}
[auto,
 block/.style ={rectangle, draw=blue, thick, fill=blue!20},
 block1/.style ={rectangle, draw=green, thick, fill=green!20},
 block2/.style ={rectangle, draw=red, thick, fill=red!20}
]
\draw (0,1)--(13,1);
\draw (0,0)--(13,0);
\path(2,0.5)node{$\dots$}; 
\path(9,0.5)node{$\dots$}; 
\path(0.5,0.5)node[block]{$1$}; 
\path(3.5,0.5)node[block]{$n-1$};
\path(6.5,0.5)node[block]{$n+1$};
\path (5,0.5)node[block]{$n$};
\path (11,0.5)node[block]{$N$};
\path (12,0.5)node{$\to$};
\end{tikzpicture}
\caption{Vehicles in  {single-lane} road.}\label{1corsia}
\end{figure}

The dynamical equations of the system are obtained combining two interaction terms.  The first one is the interaction term related to the model \eqref{bando_classico}  \cite{Bando1994203,Bando19951035}. It is a relaxation term towards a desired velocity function $V(\cdot)$ that depends only on the headway  $\Delta x_n=x_{n+1}-x_{n}>0$ between the vehicle $n$ and the vehicle ahead with index $n+1$, as shown in Fig. \ref{1corsia}. The acceleration of each vehicle is regulated by the difference between its velocity and the optimal velocity. The optimal velocity function is typically a monotonically increasing function of the headway  and it is bounded. It tends to zero for small headways and to a maximum value $V^{max}$ for large headways. Furthermore we assume that $V$ is non-negative. This term is multiplied by a parameter $\alpha_n$ denoting the speed of response of each driver, with dimensions one over time. The second term is the classical Follow-the-Leader interaction term \cite{Reuschel1950193,Reuschel2} from model \eqref{ftl_classico}, multiplied by a parameter $\beta_n$ with dimensions  {length square} over time. In this term the acceleration of a vehicle is directly proportional to the difference between the velocity of the vehicle in front and its own and is inversely proportional to their mutual distance.  

Since we are considering identical vehicles we assume $\alpha_n=\alpha$ and $\beta_n=\beta$ for all $n=1,\dots,N$.

The model is given by
\begin{equation}\label{bftl}
\begin{cases}
\dot x_n=v_n \\
\dot v_n = \alpha (V(\Delta x_n) -v_n)+\beta \dfrac{\Delta v_n}{(\Delta x_n)^2}
\end{cases}
\end{equation}
with $\Delta x_n=x_{n+1}-x_{n}$ and $\Delta v_n=v_{n+1}-v_{n}$.

In our study we usually refer to a circular road which means to solve \eqref{bftl} with periodic boundary conditions,   {in this way the vehicle with index $n=N+1$ coincides with vehicle with index $n=1$}. If we deal with a straight road we simply add an equation describing the dynamics of the leader vehicle, which must be known.

\subsection{Stability} 
Let us characterize the equilibrium for the  {single-lane} model.

\begin{proposizione}
The equilibrium of the system \eqref{bftl} is given if all vehicles are equally spaced and move with the same constant velocity.
\end{proposizione}

In fact, let us indicate with $h=\frac{L}{N}$ the constant spacing of two successive vehicles, where $L>0$ is the length of the road. Then solving \eqref{bftl} with initial  conditions 
\begin{equation}
\begin{cases}
x_{n+1}(0)-x_{n}(0)=h \\
v_n(0)=V(h) & \text{for  } n=1,\dots,N
\end{cases}
\end{equation}   {with  $x_{N+1}(\cdot)=x_1(\cdot)$ by boundary conditions},  we easily obtain the solution of the system that represents the steady state described above:
\begin{equation}\label{eq1c}
\bar{x}_n(t)=hn+V(h)t.
\end{equation} Note that the equation depends parametrically by the given number $N$ of vehicles which is constant due to the periodic boundary conditions.

Now we study the stability of model \eqref{bftl} around the equilibrium \eqref{eq1c} by linearizing the original system. Let $y_n$ be a small perturbation from the steady state \eqref{eq1c} and consider
\begin{equation}\label{pert}
x_n=\bar{x}_n+y_n.
\end{equation} Disregarding  terms higher than  $O(y_n^2)$ we obtain the linearized equation of \eqref{bftl}
\begin{equation}\label{bftl_lin}
\ddot y_n = \alpha (V'(h)\Delta y_n - \dot y_n)+\beta \frac{\Delta \dot y_n}{h^2}
\end{equation}
where $\Delta y_n=y_{n+1}-y_n$ and $\Delta \dot y_n=\dot y_{n+1}-\dot y_n$,  {again vehicle with index $n=N+1$ coincides with the vehicle with index $n=1$}.

We solve \eqref{bftl_lin} looking for solutions 
\begin{equation}
y_k(n,t)=\exp\{ia_kn+zt\}
\end{equation}
where $e^{i a_k n}$ is the Fourier coefficient with $a_k=\frac{2\pi}{N}k$, $k=0,\dots,N-1$ and $z\in \mathbb{C}$. Substituting in \eqref{bftl_lin} we obtain an equation for $z=u+iv$ 
\begin{equation}\label{zeq}
z^2+z \left(\alpha-\frac{\beta}{h^2}(e^{ia_k}-1)\right)-\alpha V'(h)(e^{ia_k}-1)=0.
\end{equation} 
If the amplitude of $y_k(n,t)$ blows up in time then the solution is unstable, so in order to find stable solutions we require that $\Re(z)=u<0$.

Let us write the two solutions of \eqref{zeq} as $z_j=u_j+iv_j$ for $j=1,2$, then the following relations holds:

$$
\begin{array}{l}
\Re(z_1+z_2)=u_1+u_2=-\alpha+ \frac{\beta}{h^2}(\cos(a_k)-1) \\
\Im(z_1+z_2)=v_1+v_2= \frac{\beta}{h^2} \sin(a_k) \\
\Re(z_1\cdot z_2)=u_1\cdot u_2-v_1 \cdot v_2=-\alpha V'(h)(\cos(a_k)-1)\\
\Im(z_1\cdot z_2)=u_1\cdot v_2-v_1 \cdot u_2= -\alpha V'(h) \sin(a_k).
\end{array}
$$
The boundary of the stability region is obtained when $u_1=0$ then $$v_1=\frac{-\alpha V'(h) \sin(a_k)}{-\alpha + \frac{\beta}{h^2}(\cos(a_k)-1)}.$$
After some algebraic manipulations we get 

\begin{equation} \label{cc}
V'(h)=\frac{\alpha}{2\cos^2\left(\frac{a_k}{2}\right)}+\frac{\beta}{h^2}+2\tan^2\left(\frac{a_k}{2}\right)\cdot \frac{\beta}{h^2}\left(\frac{\beta}{\alpha h^2}+1\right).
\end{equation}

\begin{figure}[H]
\centering
\includegraphics[scale=0.3]{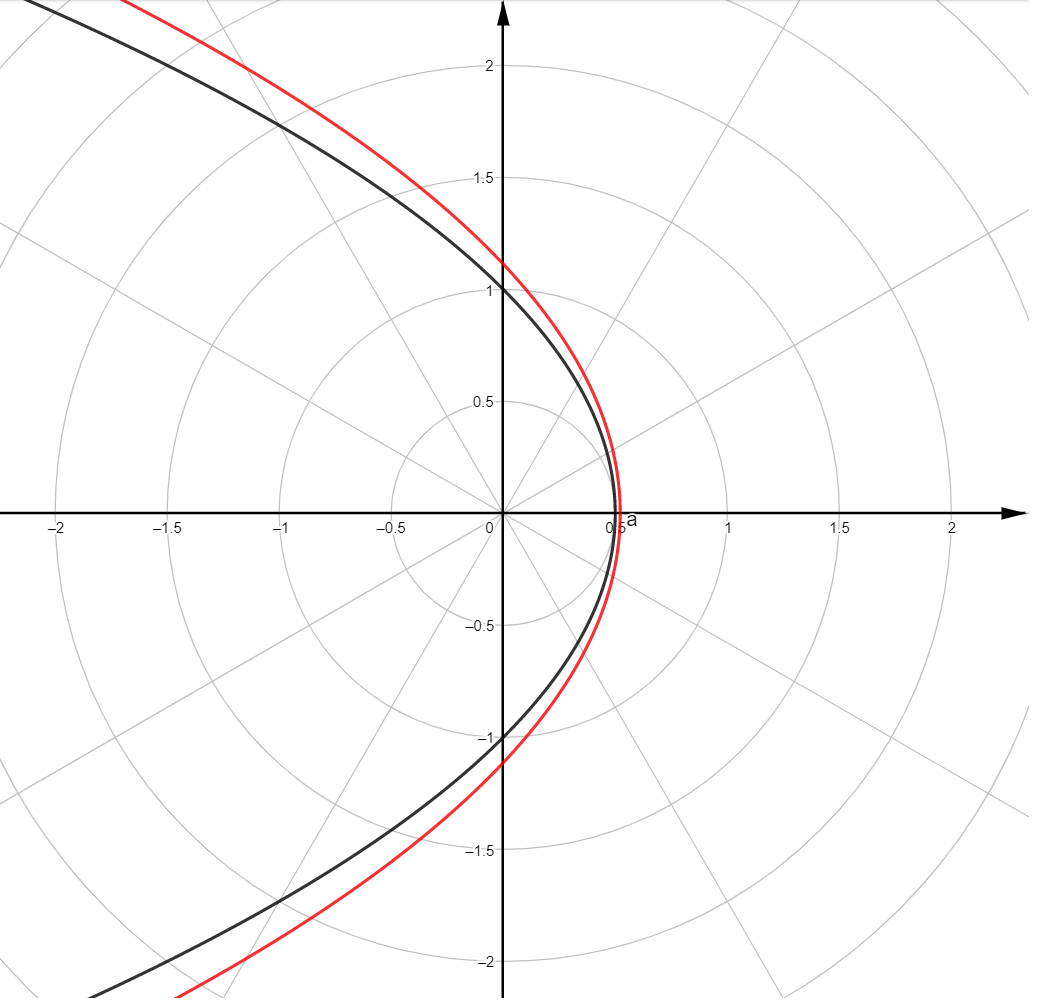}
\caption{Red: curve \eqref{cc} in the $(\alpha_k,V'(h))$ polar coordinate plane. Black: critical curve of model \eqref{bando_classico}.}
\label{critcur_plot}
\end{figure}

We can study this problem with polar coordinates in the $(\alpha_k,V'(h))$ plane as shown in Fig. \ref{critcur_plot}. The plane $(V'(h),a_k)$ can be divided into two regions: a stable region ($u<0$) and an unstable one ($u>0$) by the critical curve $u(a_k,V'(h))=0$ express by \eqref{cc}. We observe that equation \eqref{cc} coincides with the curve found in \cite{Bando19951035} if $\beta=0$. The curve \eqref{cc} is represented by the red line while the black curve is the critical curve of model \eqref{bando_classico}.


Thus we have proved the following result.

\begin{proposizione}
If \begin{equation} \label{stabilita1lane}
V'(h)<\frac{\alpha}{2}+\frac{\beta}{h^2}
\end{equation}
 {the steady state \eqref{eq1c} of model \eqref{bftl} is stable}, because for all $ k$ we have $u<0$; if $V'(h)=\frac{\alpha}{2}+\frac{\beta}{h^2}$ we have a marginal state; while  for $V'(h)>\frac{\alpha}{2}+\frac{\beta}{h^2}$ the model is unstable, because there exists at least one index $ k$ such that $u>0$.
\end{proposizione}
 For $\beta=0$ the condition \eqref{stabilita1lane}  is consistent with the stability condition derived in \cite{Bando19951035}. Remembering that $h=\frac{L}{N}$ the previous condition expresses that we gain more stability with a large number of vehicles.

\subsection{Numerical tests}
Now we present some numerical tests of model \eqref{bftl} using the Runge Kutta 5 method, with time step $\Delta t =0.1 $ s. 

\begin{figure}[H]
\centering
\includegraphics[scale=0.6]{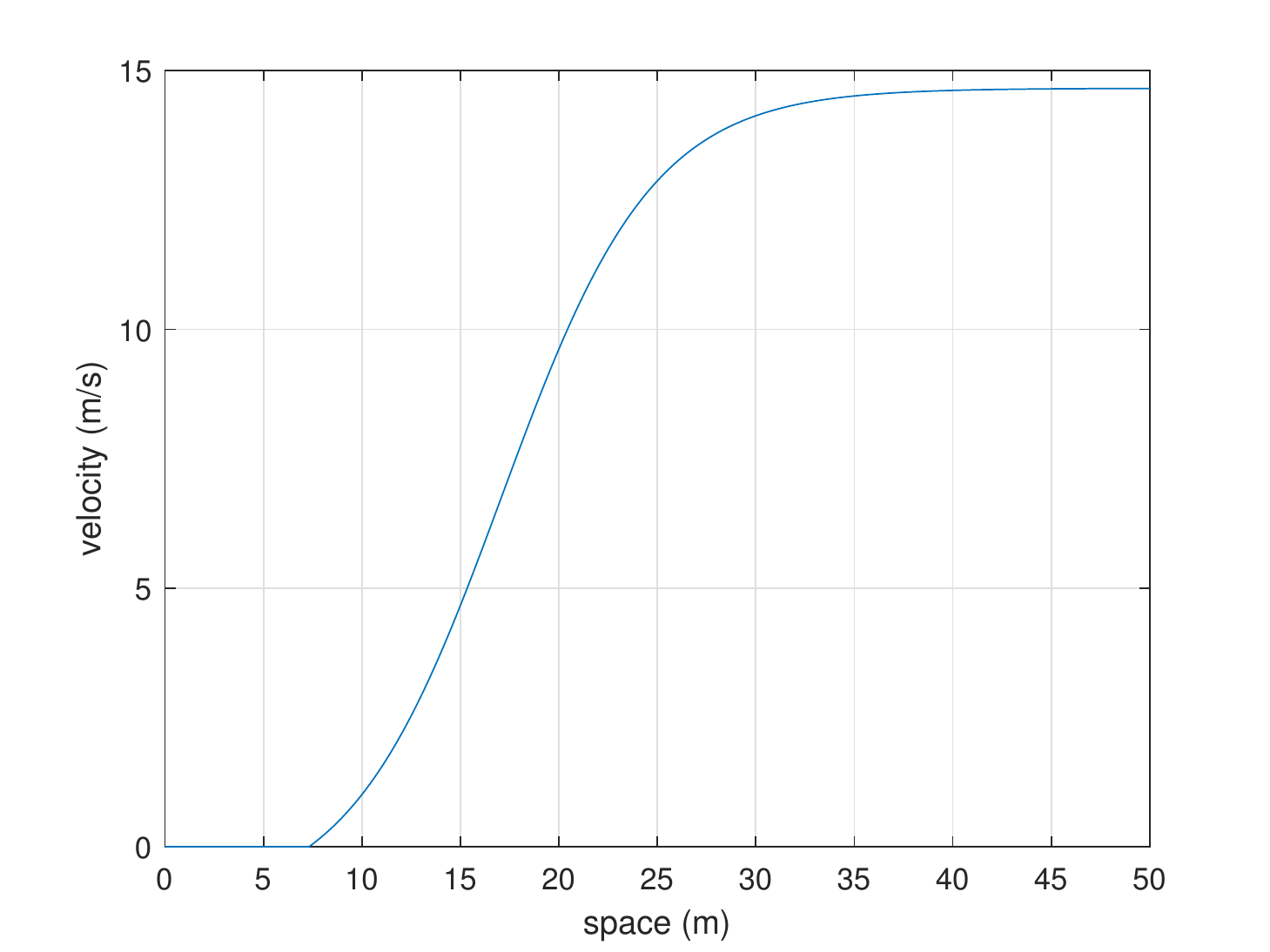}
\caption{$V(\cdot)$ function.}\label{optimalv}
\end{figure}

Let us fix $\alpha=1$ s\textsuperscript{$-1$}, $\beta=100$ m\textsuperscript{$2$}/s, $L=1500$ m,  and consider the desired velocity function expressed by \begin{equation} \label{veld}V(\Delta x) = \max\{0,V_{HT}(\Delta x)\} \end{equation} see Fig \ref{optimalv}, where \begin{equation}\label{velht} V_{HT}(\Delta x)=V_1+V_2\tanh(C_1(\Delta x -l_c)-C_2)\end{equation} is the function given by Helbing and Tilch in \cite{Helbing1998133} where they carried out a calibration of model \eqref{bando_classico} respect to the empirical data, obtaining the optimal parameter values $V_1=6.75$ m/s, $V_2=7.91$ m/s, $C_1=0.13$ m\textsuperscript{$-1$}, $C_2=1.57$ and $l_c=5$ m is the length of the vehicles. Velocity parameters $V_1,V_2$ determine the minimum expected speed $V_1-V_2$ and the maximal expected speed $V_1+V_2$, while $C_1,C_2$ are calibration parameters. Thus $V^{max}=14.66$ m/s.

\begin{figure}[H]
\centering
\includegraphics[scale=0.6]{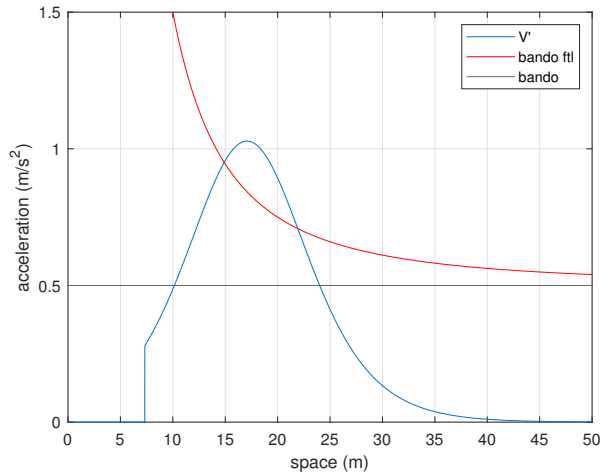}
\caption{$V'(\cdot)$ function in blue, model \eqref{bftl} stability condition in red, model \eqref{bando_classico} stability condition in black.}\label{stab_plot}
\end{figure}

From condition \eqref{stabilita1lane} we obtain that the model \eqref{bftl} with velocity \eqref{veld} is stable if $h<10.14$ m and $h>24$ m as shown in Fig.\ref{stab_plot}. In terms of number of vehicles along the circular road we have stability for $N<68$ and $N>100$. Note that, with the same parameters, the model \eqref{bando_classico} is stable for $N<62$ and $N>147$.

In the next two simulations we show a comparison between model \eqref{bftl} and model \eqref{bando_classico}, perturbing the system adding or removing a vehicle. The initial number of vehicles is chosen in such a way that the model \eqref{bftl} is stable while the model \eqref{bando_classico} is unstable according their stability condition.

\subsubsection{Test 1: adding one vehicle in the road}
In this simulation we consider $N=120$ vehicles at the equilibrium \eqref{eq1c}, equispaced with distance $\frac{L}{N}=12.5$ m and with velocities equal to $V(\frac{L}{N})$. At time $t=0$ s we perturb the system adding a one new vehicle inserting it in the position $\frac{1}{2}(x_{N}+L)$ with initial velocity equal to $V(\frac{L}{N})$. The final time is $T=1000$ s.

\medskip

Model \eqref{bftl}:
\begin{figure}[H]
\centering
\includegraphics[scale=0.5]{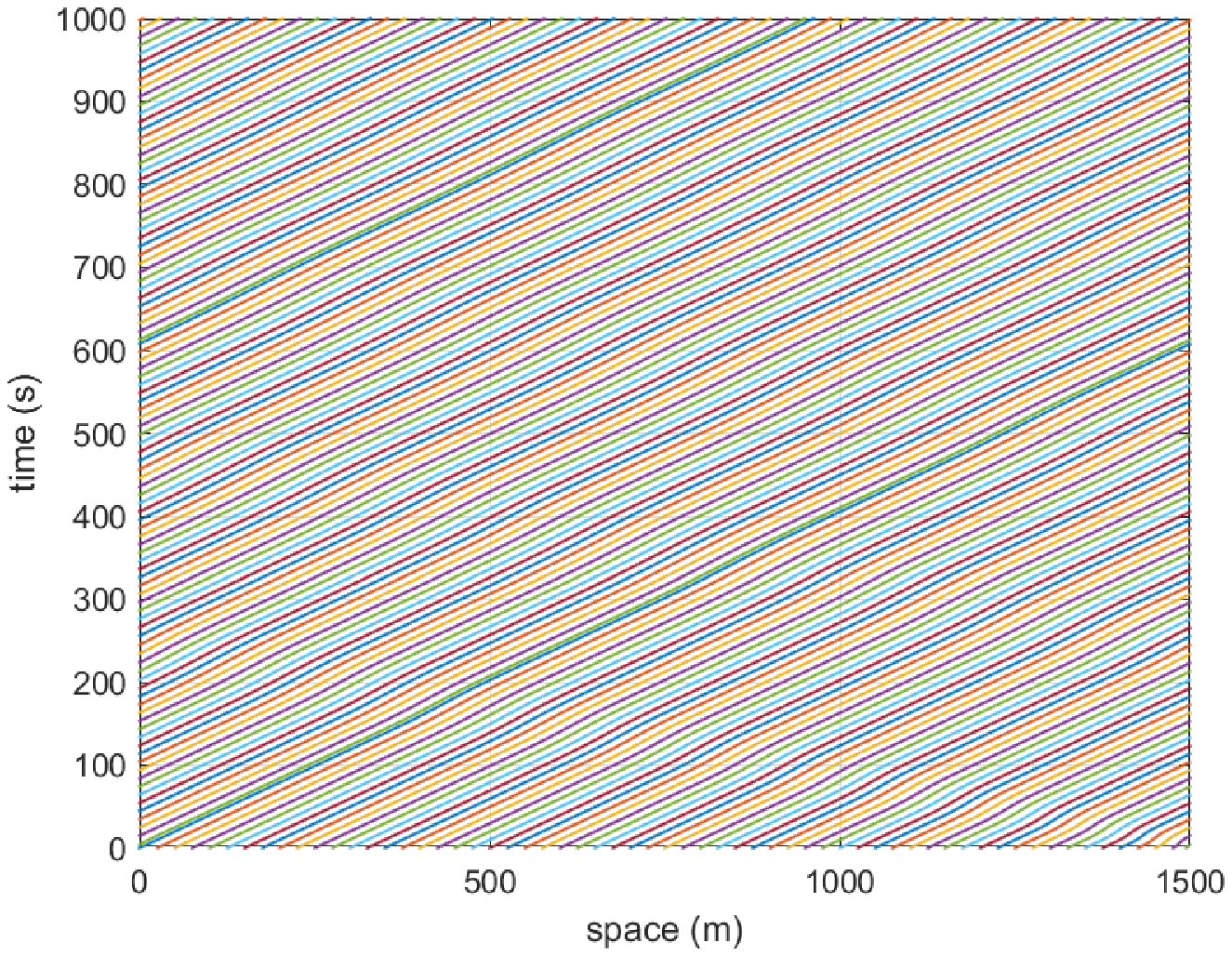} \quad \includegraphics[scale=0.5]{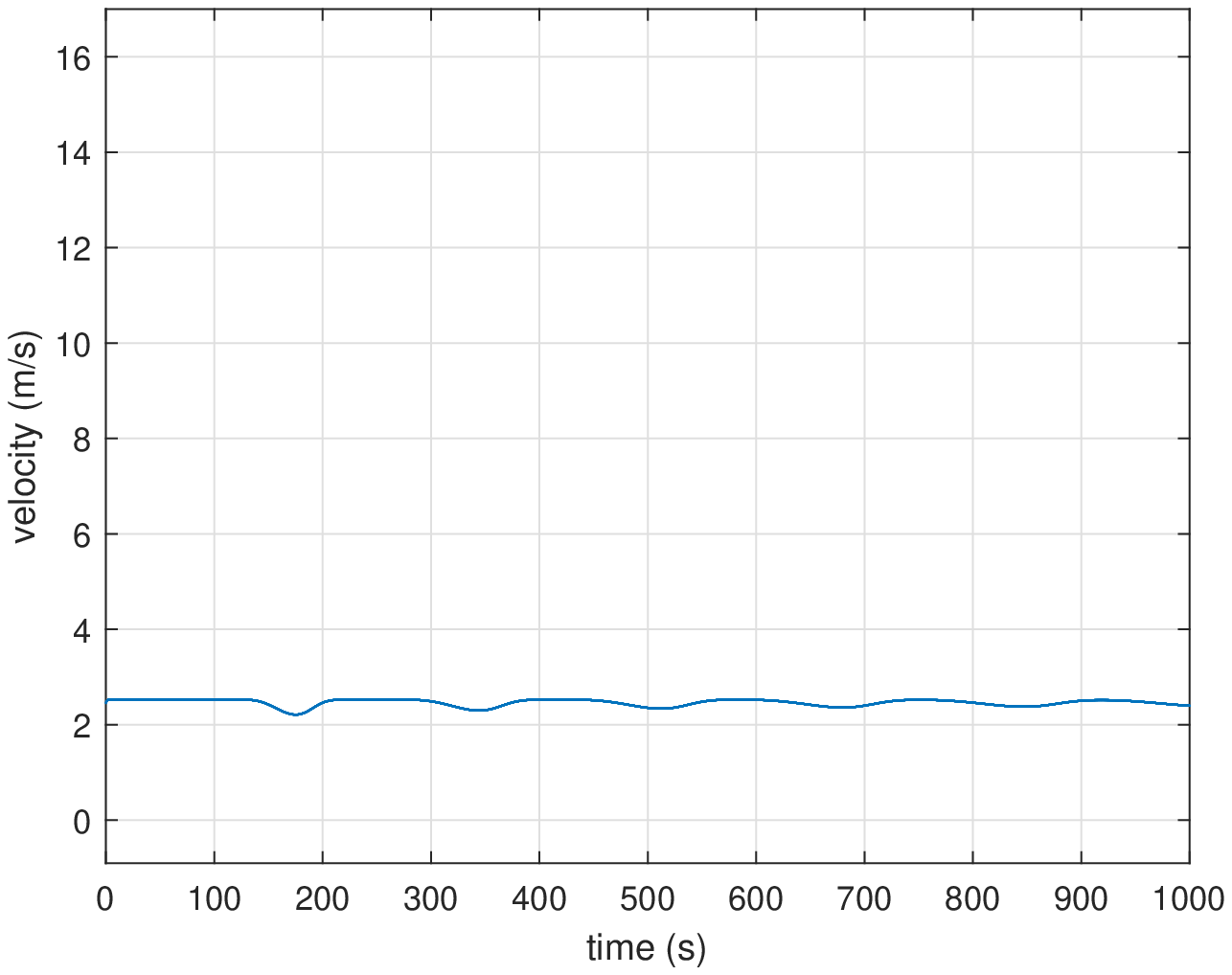}
\caption{On the left: all vehicles trajectories, on the right: velocity of vehicle 1.}
\end{figure}

Model \eqref{bando_classico}:
\begin{figure}[H]
\centering
\includegraphics[scale=0.5]{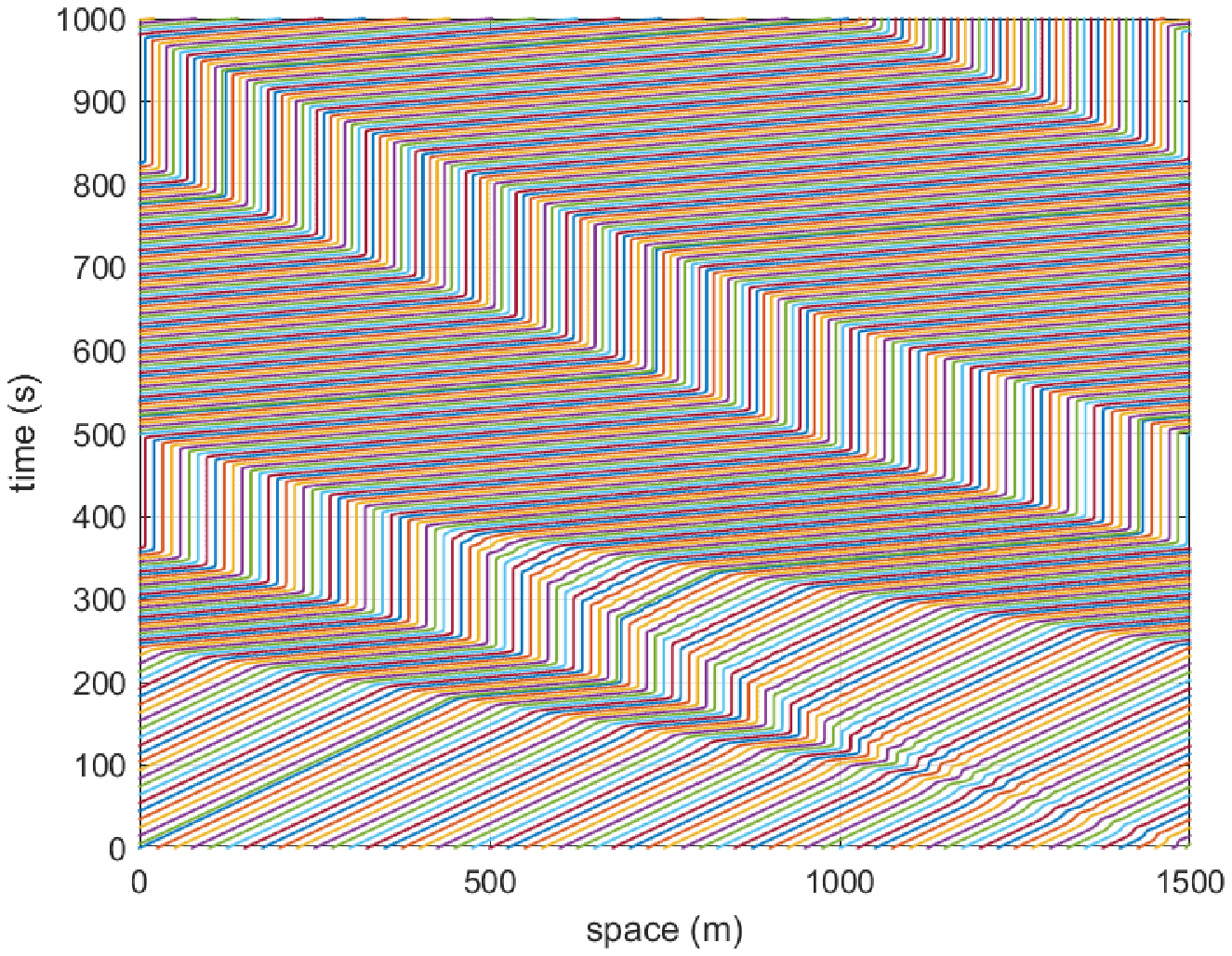} \quad \includegraphics[scale=0.5]{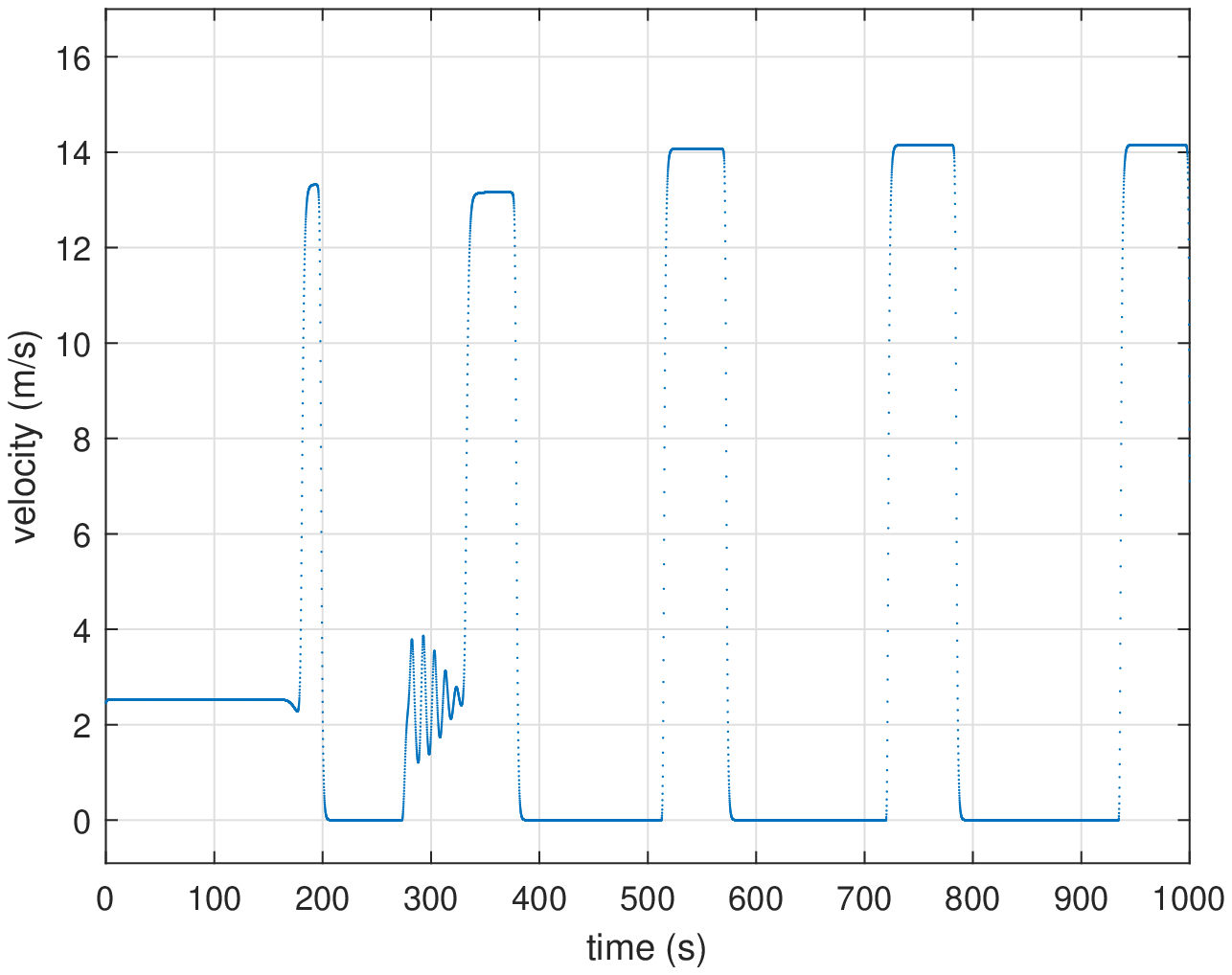}
\caption{On the left: all vehicles trajectories, on the right: velocity of vehicle 1.}
\end{figure}

We can see how the perturbation is absorbed in the in first model while it causes a creation of  {stop \& go} waves in the second model.

\subsubsection{Test 2: Removing one vehicle from the road}
In this simulation we consider again $N=120$ vehicles at the equilibrium \eqref{eq1c}, equispaced with distance $\frac{L}{N}$ and with velocities equal to $V(\frac{L}{N})$. At time $t=0$ s we perturb the system removing one vehicle choosing the one with index $N$. We set the final time $T=1000$ s.

\medskip

Model \eqref{bftl}:
\begin{figure}[H]
\centering
\includegraphics[scale=0.5]{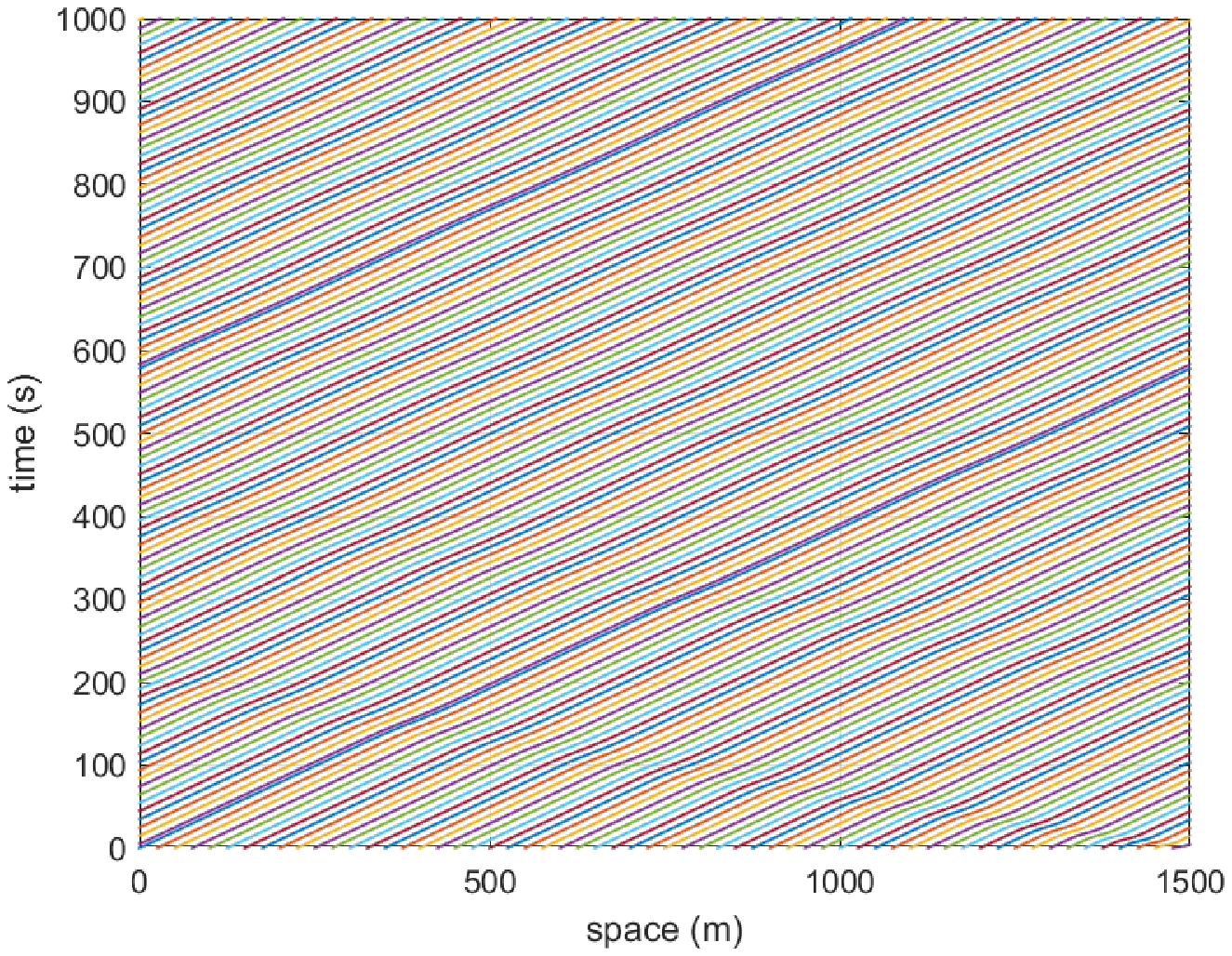} \quad \includegraphics[scale=0.5]{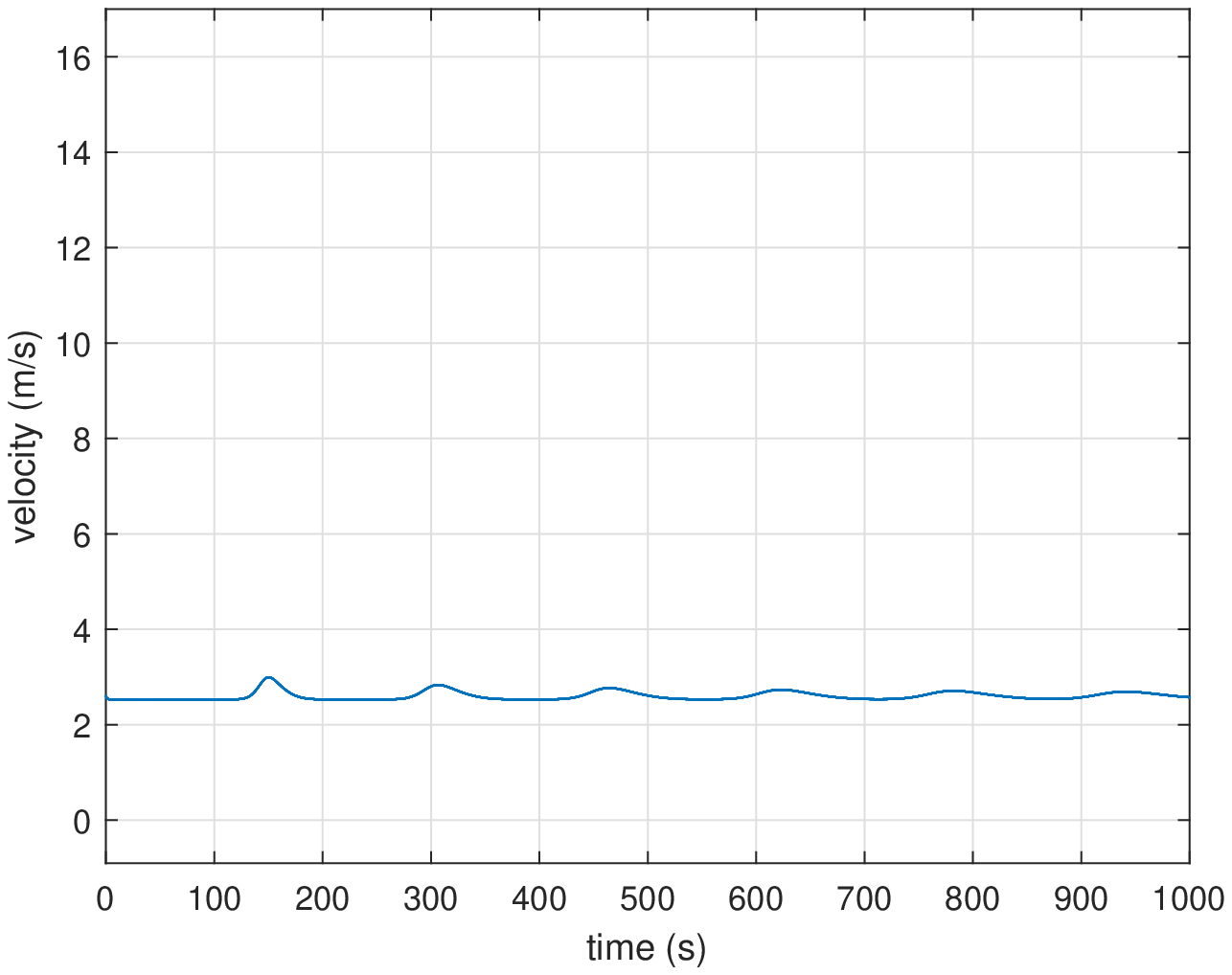}
\caption{On the left: all vehicles trajectories, on the right: velocity of vehicle 1.}
\end{figure}

Model \eqref{bando_classico}:
\begin{figure}[H]
\centering
\includegraphics[scale=0.5]{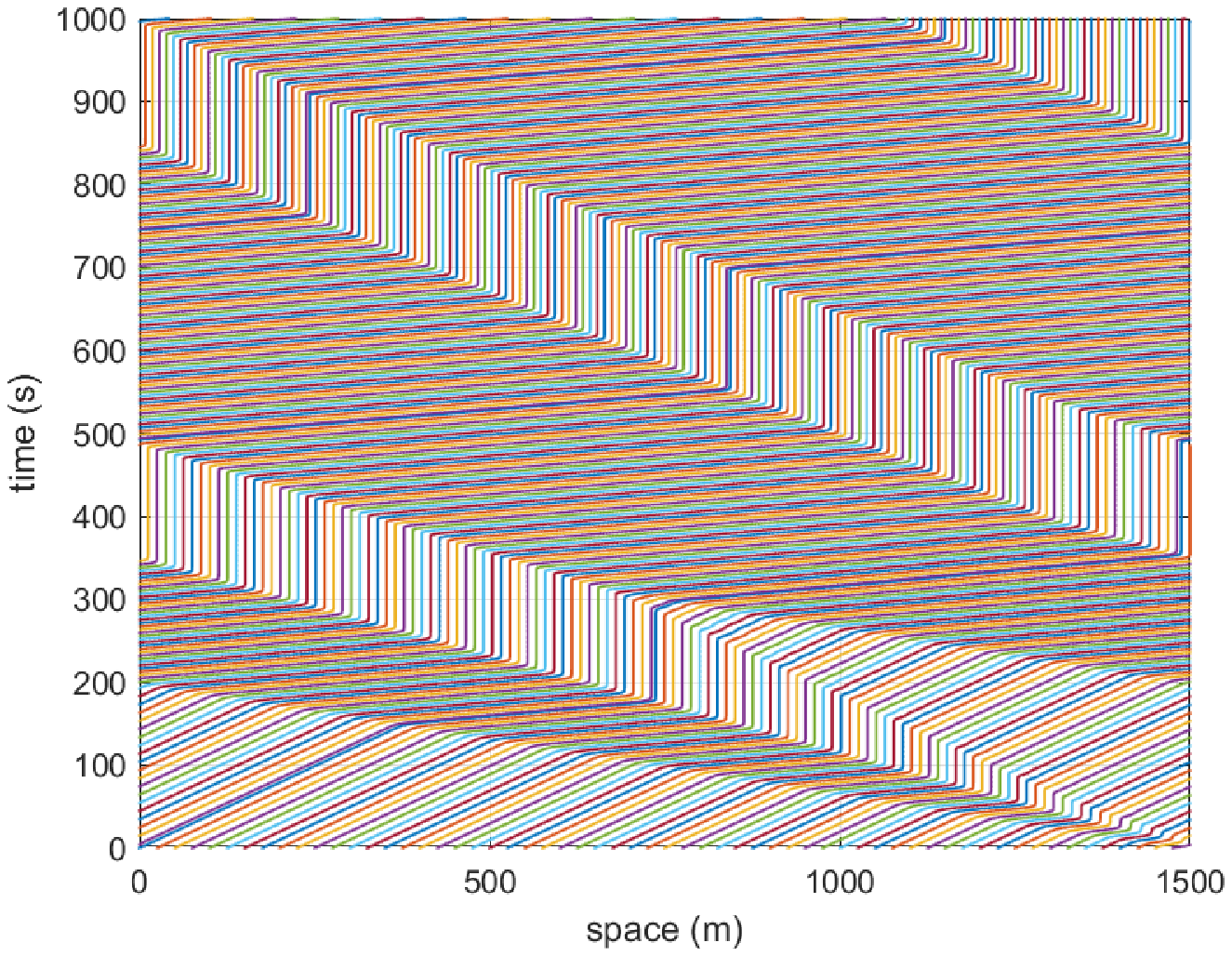} \quad \includegraphics[scale=0.5]{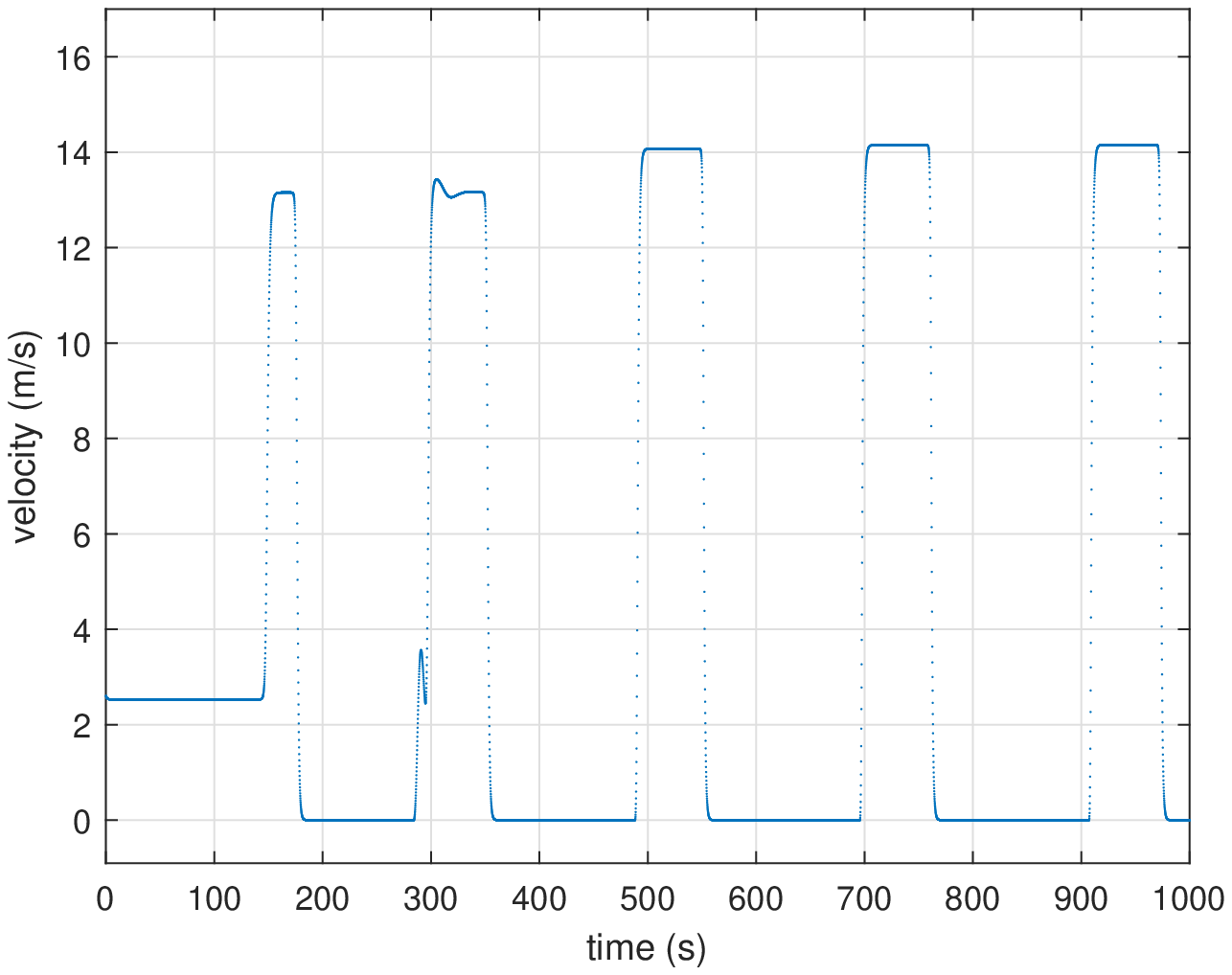}
\caption{On the left: all vehicles trajectories, on the right: velocity of vehicle 1.}
\end{figure}

Also in this test we can observe the differences when a perturbation occurs in the two models.

\subsubsection{Test 3: stop \& go waves}
In this simulation we start with $N=90$ vehicles at the equilibrium \eqref{eq1c}, equispaced with distance $\frac{L}{N}\simeq 16.66$ m and with velocities equal to $V(\frac{L}{N})$. At time $t=0$ s we perturb the system adding a new vehicle as in the previous simulations. We set $T=1000$ s.

\medskip


Model \eqref{bftl}:
\begin{figure}[H]
\centering
\includegraphics[scale=0.5]{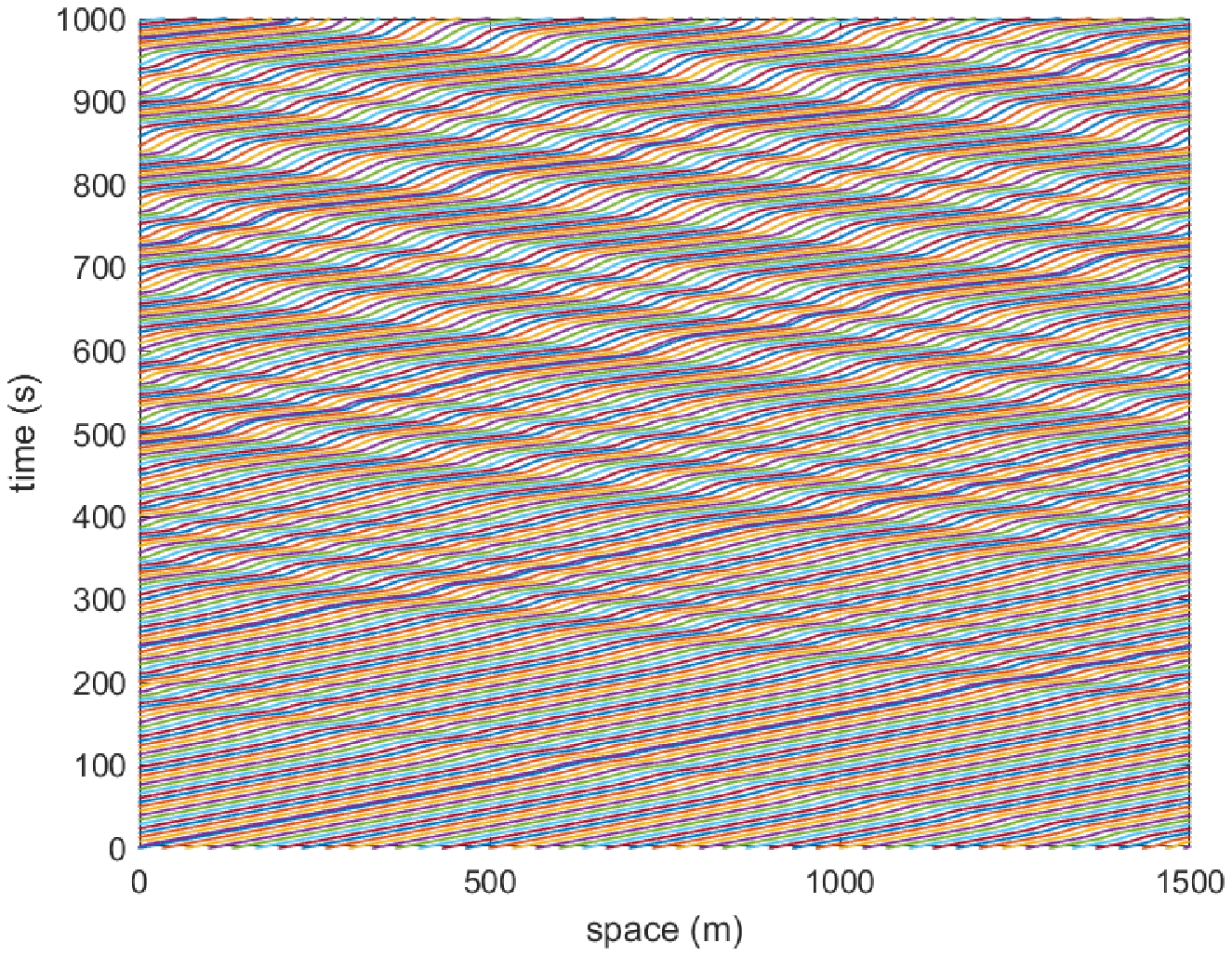} \quad \includegraphics[scale=0.5]{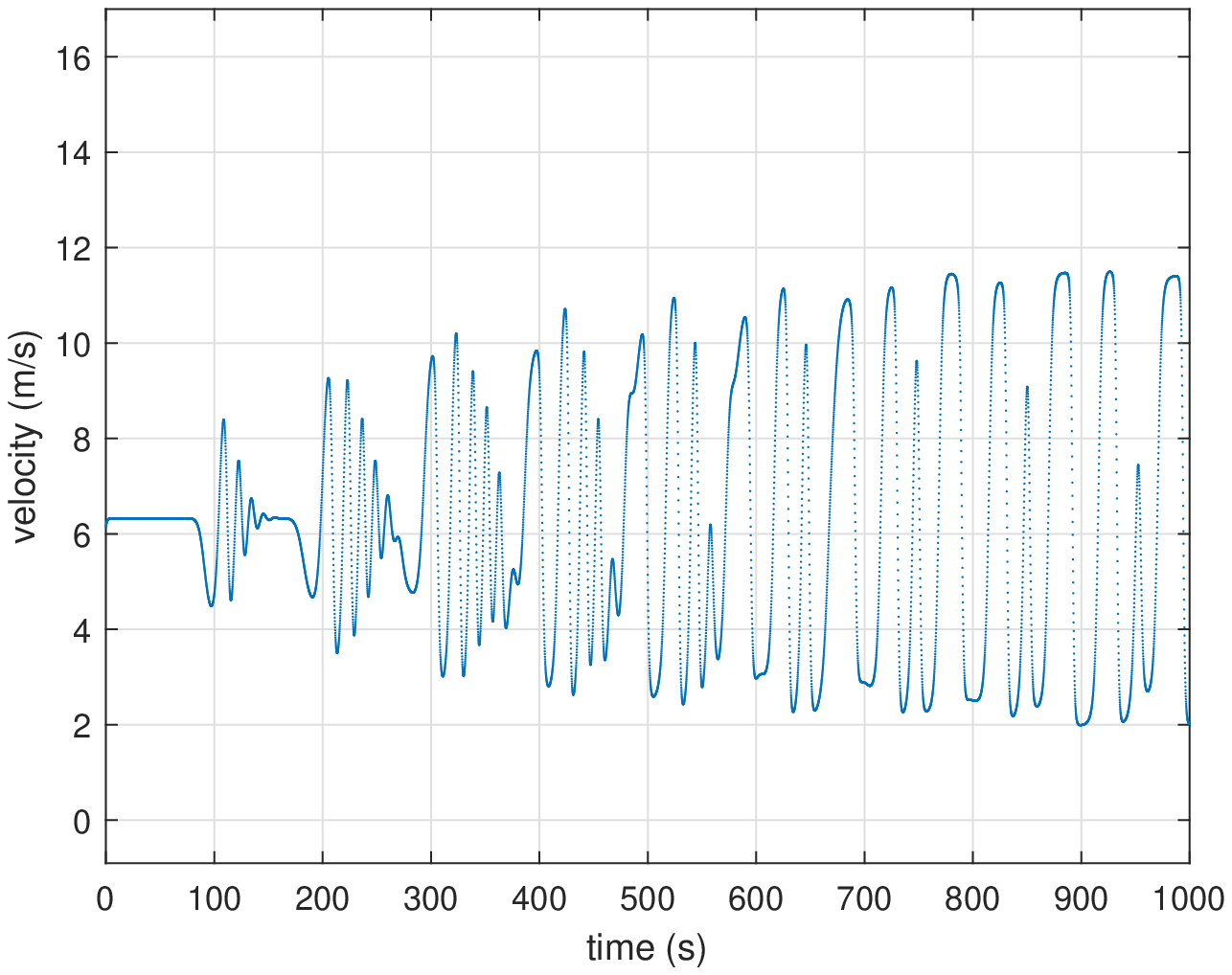}
\caption{On the left: all vehicles trajectories, on the right: velocity of vehicle 1.}
\end{figure}

Model \eqref{bando_classico}:
\begin{figure}[H]
\centering
\includegraphics[scale=0.5]{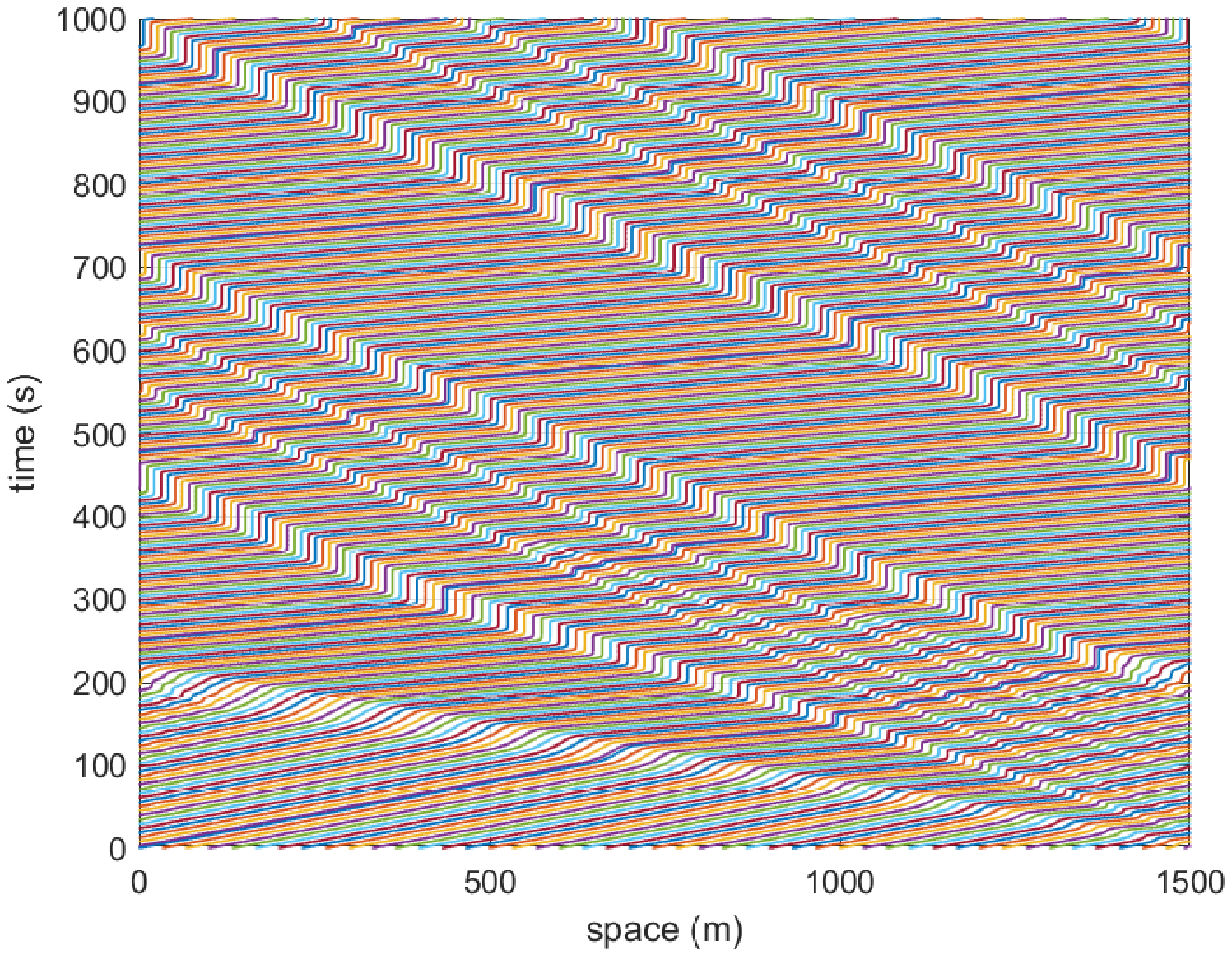} \quad \includegraphics[scale=0.5]{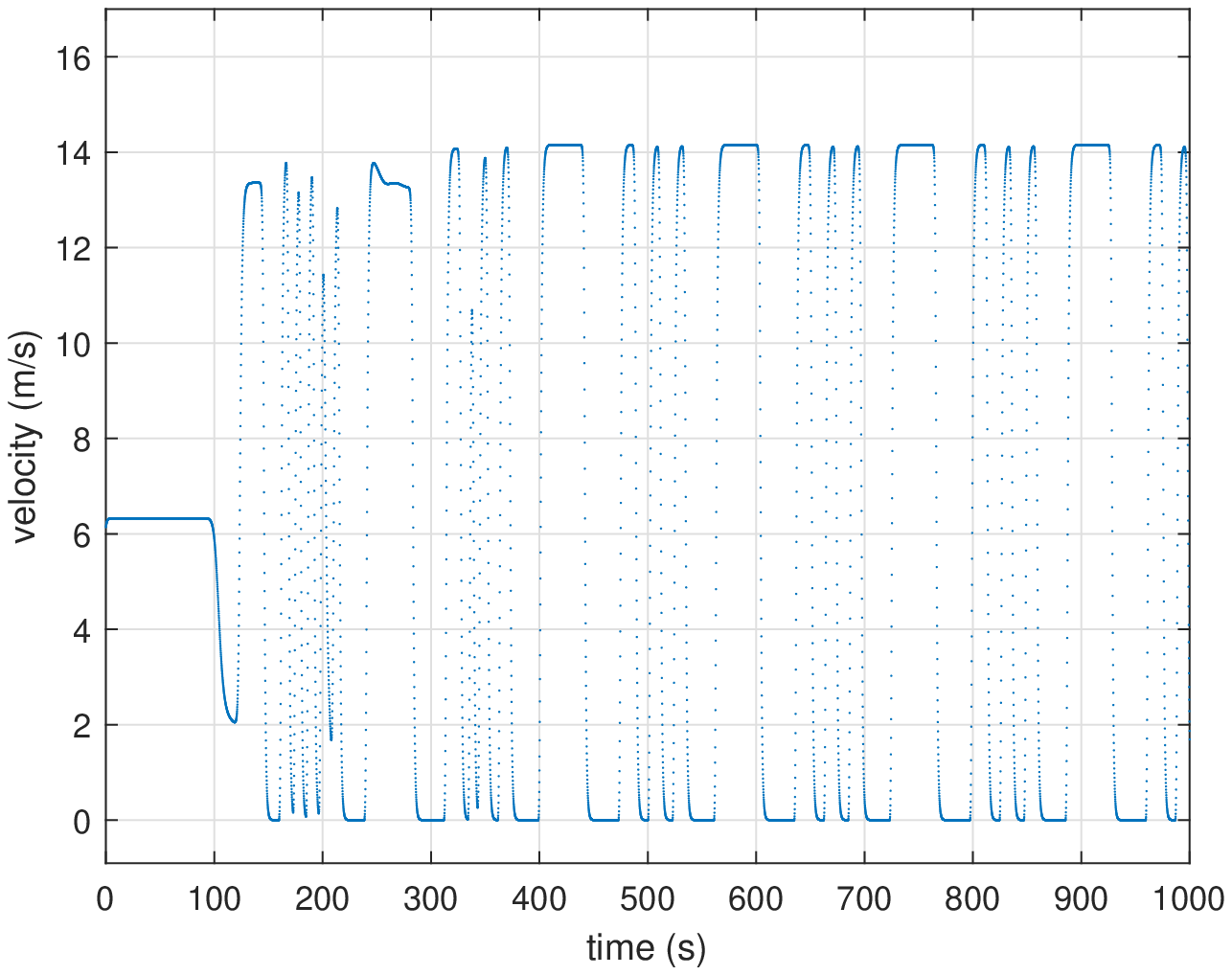}
\caption{On the left: all vehicles trajectories, on the right: velocity of vehicle 1.}\end{figure}

An example of instability for both models is reported. Although  {stop \& go} waves occur we can appreciate the differences in the oscillations of the velocity in the two models and the lack of region with zero speed in model \eqref{bftl}.

\section{Two-lane model}
\subsection{Description}
Here we study the extension of model \eqref{bftl} to a road with 2 lanes, where lane changing is allowed: lane 1 is the driving lane, while lane 2 is the fast lane. We consider a single population of homogeneous vehicles and we assume that the coefficients $\alpha,\beta$ are the same for both lanes and for all vehicles.

\begin{figure}[H]
\centering
\begin{tikzpicture}
[auto,
 block/.style ={rectangle, draw=blue, thick, fill=blue!20},
 block1/.style ={rectangle, draw=green, thick, fill=green!20},
 block2/.style ={rectangle, draw=red, thick, fill=red!20}
]
\draw (0,2)--(8,2);
\draw [dashed](0,1)--(8,1);
\draw (0,0)--(8,0);

\path(-1,0.5)node{lane 1};
\path (-1,1.5)node{lane 2};

\path(1.6,0.5)node[block]{$p_n^1$};
\path (1.2,1.5)node[block]{$p_n^2$};

\path(5.5,0.5)node[block1]{$s_n^1$};
\path (6.5,1.5)node[block1]{$s_n^2$};

\path (3.5,0.5)node[block2]{$n$};
\end{tikzpicture}
\caption{Components of the vector $\mathcal{N}_n$, containing information on cell neighbours of the $n$-th vehicle.}\label{duecorsie}
\end{figure}
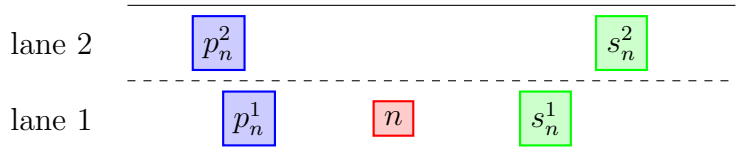

Let $N$ be the total number of vehicles in the road and $N_j=N_j(t)$ the number of vehicles in lane $j=1,2$ at time $t$; we have for all $ t,\,\, N_1(t)+N_2(t)=N $; we recall we are assuming periodic boundary conditions.  Each vehicle is identified by an index $n\in\{1,\dots,N\}$, and it is associated with a vector $\mathcal{N}_n=(j,p_n^1,p_n^2,s_n^1,s_n^2)$ whose components are: the current index lane $j\in\{1,2\}$, and the indices $s^j_n$ of the vehicle in front of vehicle $n$ in the lane $j$  (successive vehicle) and $p^j_n$  of the vehicle behind vehicle $n$ in the  lane $j$ (previous vehicle) as shown in Fig. \ref{duecorsie}. If the $n$-th vehicle does not have a successive or a previous vehicle in lane $j$ we set $s_n^j=-1$ or $p_n^j=-1$ 
respectively. In other words, the index $-1$ signifies that there is no such vehicle; for instance $s_n^1=-1$ means that the vehicle $n$ has no vehicle in front in lane 1. Whenever a lane change occurs, e.g. if the $n$-th vehicle changes lane, the vectors $\mathcal{N}_k$ for $k\in\{n,s^1_n,s^2_n,p^1_n,p^2_n\}$ affected by the change are updated with the new indices.

Assuming that vehicle $n$ is currently in lane $j$ then $\Delta x_n^j=x_{s_n^j}-x_n$ and $\Delta v_n^j=v_{s_n^j}-v_n$ denote the difference of  positions and the difference of the velocities between vehicle $n$  and its successive in the same lane. Moreover we denote with $I_j(t)=I_j$ the set of indices of vehicles ordered by their position in lane $j$ at time $t$. Note that it is sufficient to update this set only after each lane changing.

The model can be written for $j=1,2$ as

\begin{equation}\label{bftl_2lane}
\begin{cases}
\dot x_n=v_n \\
\dot v_n = \alpha (V_j(\Delta x^j_n) -v_n)+\beta \dfrac{\Delta v^j_n}{(\Delta x^j_n)^2} & n\in I_j \\

\text{ + lane changing conditions}
\end{cases}
\end{equation} where $V_j(\cdot)$ is the desired velocity function for lane $j=1,2$ with
$V^{max}_2\geq V^{max}_1$.  {In particular we assume that the velocity functions are equal to zero up to a security distance, then they monotonically increase up to their maximum value:  
\begin{equation}\begin{array}{ll}
V_1(\Delta x)=V_2(\Delta x)=0 & \Delta x \leqslant d_s \text{  (security distance)} \\
V_1(\Delta x)\leqslant V_2(\Delta x) & \text{otherwise.}
\end{array}
\end{equation}
The parameter $d_s$ is a fixed security distance that must be held by the vehicles in order to avoid collisions.}

\medskip

The lane changing rules are based according essentially on two criteria: a vehicle may change lane if it would  travel at a faster speed in the new lane, which means that is would have a higher acceleration (incentive criterion); and the changing action must be safe in order to avoid collisions with the vehicles in the adjacent lane, which means to held the security distance in every movement (security criterion).

For simplicity we introduce the compact notations: 
\begin{equation}
d(n,m)=x_m-x_n,  \qquad a_j(n,m)=\alpha (V_j(d(n,m)) - v_n)+\beta \frac{v_m-v_n}{(d(n,m))^2}
\end{equation}
to denote the difference of positions between vehicles with indices $n$ and $m$, and the acceleration of vehicle $n$ where vehicle $m$ is its successive vehicle in lane $j$.

\medskip

Thus the lane changing rules from lane $j$ to lane $j'$ can be expressed as
\begin{equation}\label{rules}
\begin{array}{ll}
a_{j'}(n,s_n^{j'})>a_j(n,s_n^j) & \text{(incentive criterion)} \\
d(n,s_n^{j'})>d_{s} \quad \text{and} \quad d(p_n^{j'},n)>d_{s} & \text{(security criterion)} 
\end{array}
\end{equation}

In particular cases we have: 
\begin{itemize}
\item if $s_n^{j'}=-1$ we consider only the security criterion;
\item if $p_n^{j'}=-1$ we consider only the incentive criterion;
\item if $s_n^{j'}=-1$ and $p_n^{j'}=-1$ we decide to change lane;
\item if $s_n^{j}=-1$  we decide to do not change lane.
\end{itemize}
Note that in this model lane changes are instantaneous and the velocity of the vehicle remains the same after the changing action. The vehicles following in the new lane adjust their velocities according to the distance from the new vehicle.

 {In order to reproduce a realistic description of traffic flow, we introduce a physical timer for lane changing because, as reported by experimental studies \cite{Kallo}, lane changing is not frequent. In other words, although a vehicle might have the opportunity and the advantage in changing lane, most often drivers prefer not to change lane. Therefore we set an expected number of lane changes per second $N_c$ and we pick randomly $N_c$ vehicles per second uniformly distributed on the set of vehicles.}

\subsection{Stability}
In the following we will use to this characterization of a steady state of model \eqref{bftl_2lane}.

\begin{proposizione}
A steady state of model \eqref{bftl_2lane} is obtained when both lanes are in equilibrium and there are no lane changing. The equilibrium velocity is given by the optimal velocity functions.
\end{proposizione}

It is easy to show that such steady state for the two-lanes model \eqref{bftl_2lane} is given when the vehicles moves with the same uniform headways $h_j=\frac{L}{N_j}$, for lane $j=1,2$ respectively, and with the optimal velocities $V_j(h_j)$. We also need to link the velocities for preserve lane changes; the condition is satisfied provided 
\begin{equation}\label{h1h2}
V_1(h_1)=V_2(h_2).
\end{equation}
Recalling that $N=N_1+N_2$, where $N$ is constant, we can write $h_2$ in terms of $h_1$ as
\begin{equation}
h_2=\frac{Lh_1}{Nh_1-L}
\end{equation} and if the equilibrium velocity is less than $V_1^{max}$ we can find a unique value for $h_1$ from equation \eqref{h1h2} that we denote by $\bar{h}_1$. Let $\bar{N_1}$ be the number of vehicles in lane 1 with headways $\bar{h}_1$ and in the same way we define $\bar{h}_2$ and $\bar{N}_2$. Thus
\begin{equation}\label{veldieq}
V^{eq}:=V_1(\bar{h}_1)=V_2(\bar{h}_2).
\end{equation} Now we prove that if \eqref{veldieq} holds we have no lane changes and both lanes remain at equilibrium. Consider model \eqref{bftl_2lane} with $\bar{N_1}$ vehicles in lane 1 and with $\bar{N}_2$ vehicles in lane 2, with initial conditions
\begin{equation}\label{equilibrium}
\begin{array}{ll}
 \forall n\in I_1 
\begin{cases}
x_n(0)\,\,\text{equally spaced with distance}\,\,\bar{h}_1 \\
v_n(0)=V^{eq}
\end{cases}\\
\forall n\in I_2
\begin{cases}
x_n(0)\,\,\text{equally spaced with distance}\,\,\bar{h}_2 \\
v_n(0)=V^{eq}.
\end{cases} 
\end{array}
\end{equation}
For the lane change from lane $1$ to lane $2$ we can show that the condition
\begin{equation}
a_{2}(n,s_n^{2})>a_1(n,s_n^1)
\end{equation}
is never verified because we have
\begin{equation}\label{conti}\begin{split}
a_{2}(n,s_n^{2})-a_1(n,s_n^1)&=\\
&=\alpha (V_2(x_{s^2_n}-x_n) - v_n)+\beta \frac{v_{s^2_n}-v_n}{(x_{s^2_n}-x_n)^2}-\alpha (V_1(x_{s^1_n}-x_n) - v_n)-\beta \frac{v_{s^1_n}-v_n}{(x_{s^1_n}-x_n)^2} \\ &=V_2(x_{s^2_n}-x_n) - V_1(x_{s^1_n}-x_n).
\end{split}
\end{equation} Moreover $\bar{h}_2<\bar{h}_1$ so the distance $x_{s^2_n}-x_n\in (d_s,\bar{h}_2-d_s)$, but from the monotonicity of the function we obtain that $V_2(h)<V_1(\bar{h}_1)\quad \forall h\in (d_s,\bar{h}_2-d_s)$. In conclusion \eqref{conti} is always negative. Similarly we can prove that there are not lane changes from lane 2 to lane 1.

We have proved the following result.

\begin{proposizione}
Consider the system \eqref{bftl_2lane} with initial conditions \eqref{veldieq}-\eqref{equilibrium}, then no lane changing occurs.
\end{proposizione}
 
\medskip 
 
In the following we study the stability of this equilibrium solution perturbing the initial headways in a lane and analysing the possibility of lane changing in both lanes. We start perturbing the slow lane (lane 1) and then the fast lane (lane 2).  Thus we start from an initial condition in which lane 1 is in a local equilibrium but does not satisfy the global equilibrium we described above. This means that we consider a uniform perturbation $\eps$ in the headways in lane 1 where we fix an initial constant headway equal to $\bar{h}_1+\eps$ and initial velocities  equal to $V_1(\bar{h}_1+\eps)$. In lane 2 we consider initial headways $\bar{h}_2$ and initial velocities  $V_2(h_2)$. We would like to study how this perturbation influences the equilibrium \eqref{equilibrium}.

\subsubsection{Case 1: perturbation in lane 1 - lane changes from lane 1 to lane 2.} We study the possibility of lane changes from lane 1 to lane 2. Let us consider a vehicle with index $n$ in lane 1, we wonder if the acceleration in lane 2 could be greater than the acceleration in lane 1
\begin{equation}
a_2(n,s_n^2)\stackrel{?}{>} a_1(n,s^1_n) \\
\Leftrightarrow
V_2(d_2)-V_1(\bar{h}_1+\eps)+\frac{\gamma}{d_2^2}(V_2(\bar{h}_2)-V_1(\bar{h_1}+\eps))\stackrel{?}{>} 0
\end{equation}
where $d_2=d(n,s^2_n)$ and $\gamma=\frac{\beta}{\alpha}$. If $\eps>0$ we do not have lane changes because the previous inequality is always false, in fact it means that in lane 1 there is now a smaller number of vehicles.
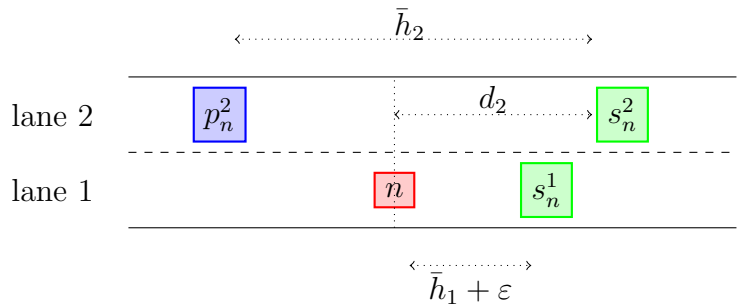
\begin{figure}[H]
\centering
\begin{tikzpicture}
[auto,
 block/.style ={rectangle, draw=blue, thick, fill=blue!20},
 block1/.style ={rectangle, draw=green, thick, fill=green!20},
 block2/.style ={rectangle, draw=red, thick, fill=red!20}
]
\draw (0,2)--(8,2);
\draw [dashed](0,1)--(8,1);
\draw (0,0)--(8,0);

\path(-1,0.5)node{lane 1};
\path (-1,1.5)node{lane 2};

\path (1.2,1.5)node[block]{$p_n^2$};

\path(5.5,0.5)node[block1]{$s_n^1$};
\path (6.5,1.5)node[block1]{$s_n^2$};

\path (3.5,0.5)node[block2]{$n$};
\draw[dotted](3.5,0)--(3.5,2);

\draw[dotted,<->](3.5,1.5)--(6.1,1.5);
\path(4.8,1.7)node{$d_2$};

\draw[dotted,<->](1.4,2.5)--(6.1,2.5);
\path(3.7,2.7)node{$\bar{h}_2$};

\draw[dotted,<->](3.7,-0.5)--(5.3,-0.5);
\path(4.5,-0.8)node{$\bar{h}_1+\eps$};

\end{tikzpicture}
\caption{Lane change from 1 to 2.}\label{cambio12}
\end{figure} 

Consider now  the case $\eps<0$, assuming $V_1$ is an invertible function, and denoted with $V_1^{-1}$ its inverse, we can write
\begin{equation}\label{eqnonap}
\eps < V_1^{-1} \left( \frac{V_2(d_2)+\frac{\gamma}{d_2^2} V_2(\bar{h}_2)}{1+\frac{\gamma}{d_2^2}} \right)-\bar{h}_1.
\end{equation}

Recalling the security criterion we have that an admissible distance $d_2$ must satisfy $d_2\in(d_s,\bar{h}_2-d_s)$ and therefore the maximum of \eqref{eqnonap} is reached when $d_2$ tends to $\bar{h}_2-d_s$. We get so this threshold for $\eps$:
\begin{equation}\label{soglianoap}
\eps < V_1^{-1} \left( \frac{V_2(\bar{h}_2-d_s)+\frac{\gamma}{(\bar{h}_2-d_s)^2} V_2(\bar{h}_2)}{1+\frac{\gamma}{(\bar{h}_2-d_s)^2}} \right)-\bar{h}_1<0.
\end{equation}

Using a Taylor expansion for $V_1$ and disregarding terms of order $O(\eps^2)$ we can also obtain an approximation at the first order of the threshold \eqref{soglianoap}. In fact the relation
\begin{equation}
V_2(d_2)-V_2(\bar{h}_2)-\eps \left(1+\frac{\gamma}{d_2^2}\right) V_1'(\bar{h}_1) \stackrel{?}{>} 0
\end{equation} is satisfied provided
\begin{equation}
 \eps<\frac{V_2(d_2)-V_2(\bar{h}_2)}{\left(1+\frac{\gamma}{d_2^2}\right) V_1'(\bar{h}_1)}.
\end{equation}
Then using the monotonicity of the velocity function we get this a priori bound, approximated at the first order respect to $\eps$ 
\begin{equation}\label{soglia1}
 \eps<\frac{V_2(\bar{h}_2-d_s)-V_2(\bar{h}_2)}{\left(1+\frac{\gamma}{(\bar{h}_2-d_s)^2}\right) V_1'(\bar{h}_1)}<0.
\end{equation}
So if $\eps$ is smaller than this value we have lane changes from lane 1 to lane 2.

\subsubsection{Case 2: perturbation in lane 1 -  lane changes from lane 2 to lane 1.} Consider a vehicle with index $n$ in lane 2 as in Fig. \ref{cambio21}. This vehicle will change to lane 1 if the following condition is satisfied
\begin{equation}
\begin{CD}
a_1(n,s_n^1)\stackrel{?}{>} a_2(n,s^2_n).
\end{CD} 
\end{equation}
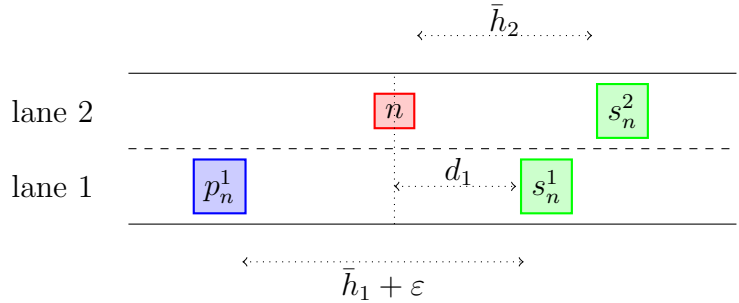
\begin{figure}[H]
\centering
\begin{tikzpicture}
[auto,
 block/.style ={rectangle, draw=blue, thick, fill=blue!20},
 block1/.style ={rectangle, draw=green, thick, fill=green!20},
 block2/.style ={rectangle, draw=red, thick, fill=red!20}
]
\draw (0,2)--(8,2);
\draw [dashed](0,1)--(8,1);
\draw (0,0)--(8,0);

\path(-1,0.5)node{lane 1};
\path (-1,1.5)node{lane 2};

\path (1.2,0.5)node[block]{$p_n^1$};

\path(5.5,0.5)node[block1]{$s_n^1$};
\path (6.5,1.5)node[block1]{$s_n^2$};

\path (3.5,1.5)node[block2]{$n$};
\draw[dotted](3.5,0)--(3.5,2);

\draw[dotted,<->](3.5,0.5)--(5.1,0.5);
\path(4.35,0.7)node{$d_1$};

\draw[dotted,<->](3.8,2.5)--(6.1,2.5);
\path(4.95,2.7)node{$\bar{h}_2$};

\draw[dotted,<->](1.5,-0.5)--(5.2,-0.5);
\path(3.35,-0.8)node{$\bar{h}_1+\eps$};

\end{tikzpicture}
\caption{Lane change from 2 to 1.}\label{cambio21}
\end{figure}
In this case clearly we will not have lane changes if $\eps<0$. Thus we consider only the case $\eps>0$ and we obtain
\begin{equation}
V_1(d_1)-V_1(\bar{h}_1)+\frac{\gamma}{d_1^2} V_1'(\bar{h}_1)\eps \stackrel{?}{>} 0
\end{equation}
where $d_1=d(n,s^1_n)$ with admissible distance $d_1\in(d_s,\bar{h}_1+\eps-d_s)$.  If $d_1>\bar{h}_1$ the previous relation is always verified, while if $d_1\leq \bar{h}_1$ considering the security criterion we can conclude that the perturbation must be greater than the safety distance in order to activate lane changes:
\begin{equation}\label{soglia22}
\eps>d_s.
\end{equation}
In fact the arrival of a vehicle from lane 2 modifies the initial perturbation $\eps$, decreasing the headway in lane 1, and we go back to the case 1.

%

\medskip

Now we repeat the same analysis adding a perturbation $\eps$ in the initial headways in lane 2 starting from the equilibrium \eqref{equilibrium}. Thus we consider an initial condition where vehicles in lane 1 have initial headways $\bar{h}_1$ and initial velocities  $V_1(\bar{h}_1)$, and vehicles in lane 2 have initial headways $\bar{h}_2+\eps$ and initial velocities  $V_2(\bar{h}_2+\eps)$.

\subsubsection{Case 3: perturbation in lane 2 - lane changes from lane 1 to lane 2.} In this case a vehicle in lane 1 could clearly have a greater acceleration from lane 1 to the lane perturbed if $\eps>0$, but from the security criterion the perturbation must be satisfy the condition
\begin{equation}
\eps>d_s
\end{equation}
as seen in case 2.

\subsubsection{Case 4: perturbation in lane 2 - lane changes from lane 2 to lane 1.} If the perturbation $\eps$ is positive we expect no lane changes of this type. Therefore let us consider the case $\eps<0$. Let $n$ be the index of a vehicle in lane 1 we wonder if
\begin{equation}
a_1(n,s_n^1)\stackrel{?}{>} a_1(n,s^2_n) \\
\Leftrightarrow
V_1(d_1)-V_2(\bar{h}_2+\eps)+\frac{\gamma}{d_1^2}(V_1(\bar{h}_1)-V_2(\bar{h}_2)+\eps))\stackrel{?}{>} 0
\end{equation}
with admissible distance $d_1\in(d_s,\bar{h}_1-d_s)$. Consider the maximum distance $d_1=\bar{h}_1-d_s$, the previous inequality is satisfy if 
\begin{equation}\label{sogliaex2}
\eps < V_2^{-1} \left( \frac{V_1(d_1)+\frac{\gamma}{(d_1)^2} V_1(\bar{h}_1)}{1+\frac{\gamma}{(d_1)^2}} \right)-\bar{h}_2<0
\end{equation}
which can be linear approximated by
\begin{equation}
 \eps<\frac{V_1(d_1)-V_1(\bar{h}_1)}{\left(1+\frac{\gamma}{d_1^2}\right) V_2'(\bar{h}_2)}.
\end{equation}
Then using the monotonicity of the velocity function we get this a priori bound, approximated at the first order respect to $\eps$ 
\begin{equation}\label{soglia2apo}
 \eps<\frac{V_1(\bar{h}_1-d_s)-V_1(\bar{h}_1)}{\left(1+\frac{\gamma}{(\bar{h}_1-d_s)^2}\right) V_2'(\bar{h}_2)}<0.
\end{equation}

\medskip
We can summarize the results in the following proposition.

\begin{proposizione}
Starting from the equilibrium, lane changing for system \eqref{bftl_2lane} are activated if a perturbation $\eps$ in the headways satisfies the thresholds in Tab.~\ref{tabella1}. Therefore there are perturbations that do not affect the equilibrium of the system.
\end{proposizione}

\begin{table}[h]
\centering
\begin{tabular}{ |p{5cm}|c|c|  }
\hline
 & from lane 1 to lane 2 & from lane 2 to lane 1 \\
\hline 
\rule{0mm}{1cm}
 perturbation $\eps$ in lane 1 (slow lane) & $ \eps<\dfrac{V_2(\bar{h}_2-d_s)-V_2(\bar{h}_2)}{\left(1+\frac{\gamma}{(\bar{h}_2-d_s)^2}\right) V_1'(\bar{h}_1)}<0$ & $\eps>d_s>0$\\
 
\hline
\rule{0mm}{1cm}
perturbation $\eps$ in lane 2 (fast lane) & $\eps>d_s>0$ &  $\eps<\dfrac{V_1(\bar{h}_1-d_s)-V_1(\bar{h}_1)}{\left(1+\frac{\gamma}{(\bar{h}_1-d_s)^2}\right) V_2'(\bar{h}_2)}<0$\\

\hline
\end{tabular}
\caption{\label{tabella1}Thresholds and perturbations.}
\end{table}

\subsection{Numerical tests}
Here we present some numerical tests for the two-lane model \eqref{bftl_2lane}, using the Runge Kutta 5 method. In the following simulations we set a maximum number of lane changes per second equal to $N_c=1$ and we fix $\Delta t=0.1$ s.

Let us set $L=1500$ m, $\alpha=5$ s\textsuperscript{$-1$}, $\beta=100$ m\textsuperscript{2}/s. We use the two optimal velocity functions defined in   \eqref{veld} with parameters $V_1=0$, $V_2=5$, $C_1=0.02$ m\textsuperscript{$-1$}, $C_2=0$, $l_c=5$ m, thus
\begin{equation}\label{vel2corsie}
V_1(h)= \begin{cases} 
5\tanh (0.02(h-5)) & \text{if   } h>d_s \\
0 & \text{otherwise}
\end{cases} 
\qquad
V_2(h)=2V_1(h).
\end{equation}
with $d_s=5$ m. We make this choice in order to verify the stability condition \eqref{stabilita1lane} in both single lanes for every value of $N$. We are interesting to study the stability of the model due to the lane changes. 

\begin{figure}[H]\label{2vel}
\centering
\includegraphics[scale=0.5]{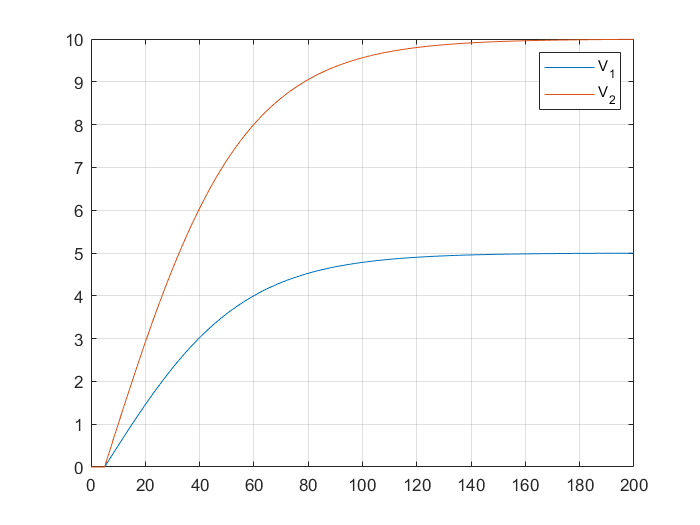}
\caption{Optimal velocity functions.}
\end{figure}

\subsubsection{Test 1: perturbation in lane 1 and lane changing from lane 1 to lane 2}
In this simulation we want to study the perturbation of the lane one from the equilibrium state. Let us fix $N=100$. Solving equation \eqref{h1h2} we get the values $\bar{h}_1=45.4$ m and $\bar{h}_2=22.4$ m for which the system remains at the equilibrium if we start from the corresponding steady state. In this case $\bar{N}_1=33$ m and $\bar{N}_2=67$ m.

Now we want to perturb the lane 1 adding new vehicles. From bound \eqref{soglia1} we obtain that the perturbation $\eps$ in the headways of lane 1 that enables lane changing from lane 1 to lane 2 must satisfy $\eps<-16.5$ m, which means that lane changes occur only if $N_1>51.7$.

Thus fix $\tilde{\eps}=-16.59$ m in order to have $N_1(0)=52$ and set $N_2(0)=\bar{N}_2$. We consider the following initial data
\begin{equation}
\begin{array}{ll}
 \forall n\in I_1 
\begin{cases}
x_n(0)\,\,\text{equally spaced with distance}\,\,\bar{h}_1+\tilde{\eps} \\
v_n(0)=V_1(\bar{h}_1+\tilde{\eps})
\end{cases}\\
\forall n\in I_2
\begin{cases}
x_n(0)\,\,\text{equally spaced with distance}\,\,\bar{h}_2 \\
v_n(0)=V_1(\bar{h}_2)
\end{cases} 
\end{array}
\end{equation}

\begin{figure}[H]
\centering
\includegraphics[scale=0.4]{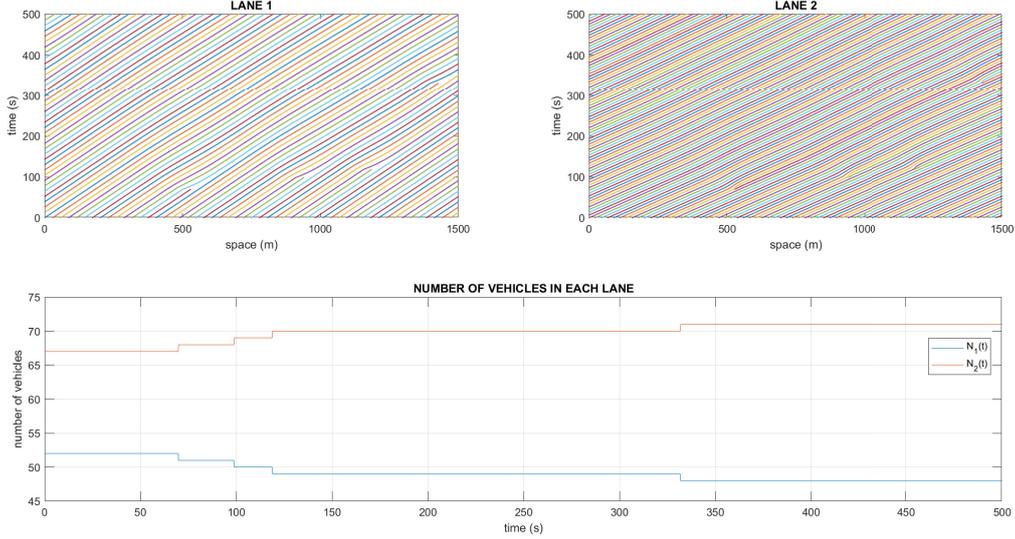}
\caption{Top: vehicle trajectories in the two lanes. Bottom: number of vehicles versus time.}\label{cambicorsia}
\end{figure}

Fig. \ref{cambicorsia} shows the simulation for $T=500$ s. We can see how the perturbation in lane 1 causes lane changes to lane 2 as expected, until the number of vehicles in lane 1 is such that the headways become smaller than the value $\bar{h}_1+\eps$ for which we cannot have any more lane changes. In this particular case a new equilibrium is reached with $N_1(T)=48$ and the corresponding headways in lane 1 are equal to $\frac{L}{N_1(T)}=31.25=\bar{h}_1-13.48$ m. This corresponds to a perturbation with $\eps=-14.15$ m which is greater than the threshold above. Thus no more lane changes are expected and the system has acquired a new equilibrium with $h_1=31.25$ m $h_2=21.13$ m and $V_1(h_1)=2.41$ m/s, $V_2(h_2)=3.12$ m/s.
 
\subsubsection{Test 2: perturbation in lane 1 and lane changing from lane 2 to lane 1} 
Whit this simulation we want to study the possibility of lane changes from lane 2 to lane 1. We consider again the equilibrium found in Test 1, and we focus attention to perturb the headways in lane 1 with a positive value of $\eps$, which means to remove some vehicles from the initial value $\bar{N}_1$.

From \eqref{soglia22} we know that a perturbation that activates lane changes from lane 2 to lane 1 must be greater that the security distance. In our case this is verify if we consider $N_1<29.73$ vehicles at initial time. Therefore we fix $\tilde{\eps}=6.27$ m in order to have $N_1(0)=29$ and set $N_2(0)=\bar{N}_2$. Thus the initial conditions are given by
\begin{equation}
\begin{array}{ll}
 \forall n\in I_1 
\begin{cases}
x_n(0)\,\,\text{equally spaced with distance}\,\,\bar{h}_1+\tilde{\eps} \\
v_n(0)=V_1(\bar{h}_1+\tilde{\eps})
\end{cases}\\
\forall n\in I_2
\begin{cases}
x_n(0)\,\,\text{equally spaced with distance}\,\,\bar{h}_2 \\
v_n(0)=V_1(\bar{h}_2)
\end{cases} 
\end{array}
\end{equation}
 
\begin{figure}[H]
\centering
\includegraphics[scale=0.4]{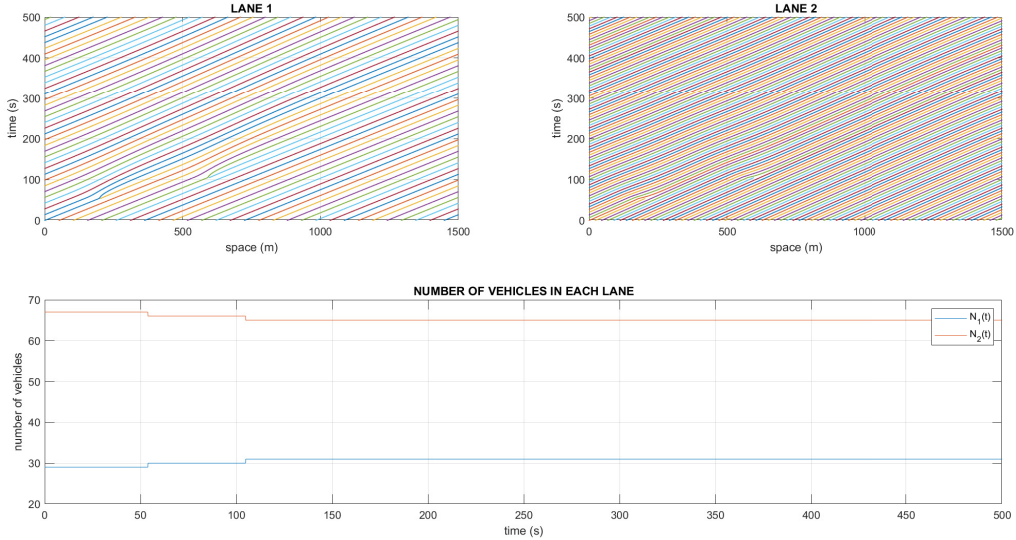}
\caption{Top: vehicle trajectories in the two lanes. Bottom: number of vehicles versus time.}\label{cambicorsia2}
\end{figure} 
 
Fig. \ref{cambicorsia2} shows the simulation for $T=500$ s. We can see how the perturbation in lane 1 causes lane changes from lane 2 to lane 1 as predicted. A new equilibrium is reached with $N_1(T)=31$ and the corresponding headways in lane 1 are equal to $\frac{L}{N_1(T)}=48.39=\bar{h}_1+2.99$ m. This corresponds to a perturbation with $\eps=2.99$ m which is smaller than the threshold above. Thus no more lane changes are expected and the system has acquired a new equilibrium with $h_1=48.38$ m $h_2=23.07$ m and $V_1(h_1)=3.50$ m/s, $V_2(h_2)=3.46$ m/s. 
 
\subsubsection{Test 3: evolution towards equilibrium}
In this simulation we study the evolution towards equilibrium. We start with the same number of vehicles in both lanes $N_1(0)=N_2(0)=50$. At the initial time all vehicles are equally spaced with zero velocity. 

\begin{figure}[H]
\centering
\includegraphics[scale=0.4]{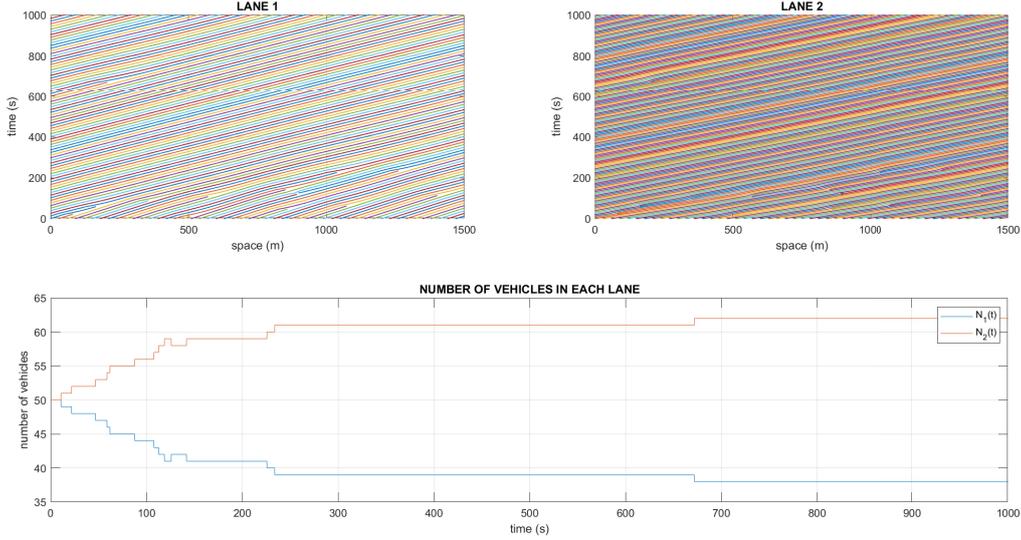}
\caption{Top: vehicle trajectories in the two lanes. Bottom: number of vehicles versus time.}\label{simu2}
\end{figure}

Fig. \ref{simu2} shows the simulation for $T=1000$ s. We can see the presence of an initial phase where vehicles change lane more frequently until arriving in a phase with few lane changes that let the traffic more regular. Initially all vehicles accelerate and lane 1 is partially defected by lane changes towards lane 2 until $N_1(T)=38,N_2(T)=62$. In this simulation lane changes from lane 1 to lane 2 are the 92.8\% of the total lanes changes. 

\subsubsection{Test 4: stop \& go waves\label{testsg}}
In this simulation we use the two velocity functions as in \eqref{veld} in order to consider also the instability due to the number of vehicles as seen in the  {single-lane} case. We fix 
$\alpha=1,\beta=100$ and $$ V_1(\Delta x)=
\begin{cases}
6.75+7.91\tanh(0.13(\Delta x -5)-1.57) & \Delta x >5 \\
0 & \text{otherwise}
\end{cases} \qquad V_2(\Delta x)=2V_1(\Delta x).
$$ 
The stability conditions for the  {single-lane} \eqref{stabilita1lane} are in this case: for lane 1 stability for $N<68$ and $N>100$, while for lane 2 we have stability for $N<57$ and $N>130$.

We start with the same number of vehicles in both lanes $N_1(0)=N_2(0)=90$; lane 2 is at the equilibrium while in lane 1 we add random perturbations $r_n$ in the initial positions of the vehicles. Thus we have
\begin{equation}
\begin{array}{ll}
 \forall n\in I_1 
\begin{cases}
x_n(0)-x_{n-1}(0)=\frac{L}{N_1(0)}+r_n \\
v_n(0)=V_1(\frac{L}{N_1(0)})
\end{cases}\\
\forall n\in I_2
\begin{cases}
x_n(0)-x_{n-1}(0)=\frac{L}{N_2(0)} \\
v_n(0)=V_2(\frac{L}{N_2(0)})
\end{cases} 
\end{array}
\end{equation}

\begin{figure}[H]
\centering
\includegraphics[scale=0.4]{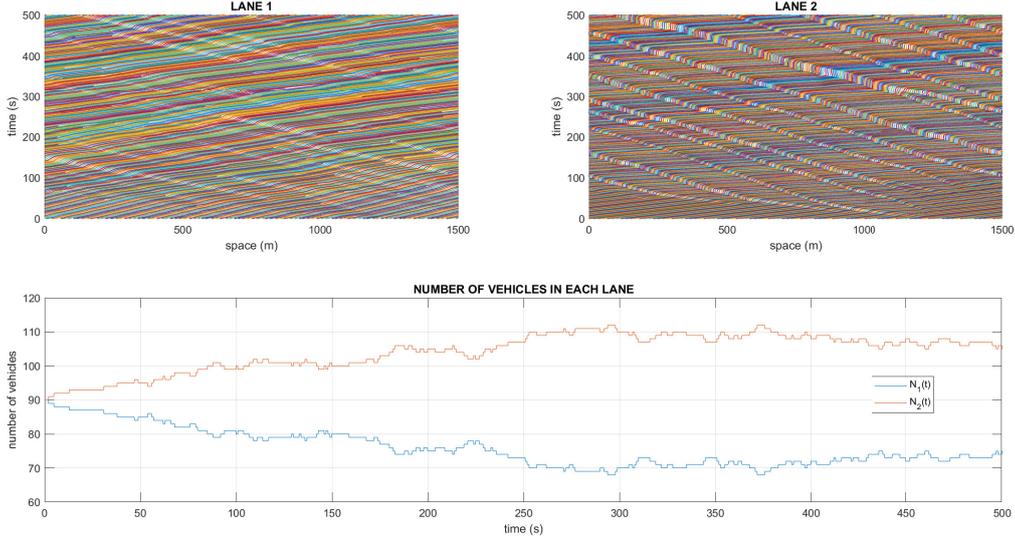}
\caption{Top: vehicle trajectories in the two lanes. Bottom: number of vehicles versus time.}\label{simi5}
\end{figure}

Fig. \ref{simi5} shows the simulation for $T=500$ s. We can see the creation of  {stop \& go} waves in both lanes due to the frequently lane changes and to the instability of the model. 

%

\section{Generalization to the multi-lane case}
The model \eqref{bftl_2lane} can be easily generalized to the  {multi-lane} case with a generic number of lanes. We can differentiate the lanes by attributing different profiles of desired velocity, therefore let $J$ be the number of lanes, we consider the velocities functions $V_1(\cdot),\dots,V_J(\cdot)$ with the property $V_i(\cdot)\leqslant V_j(\cdot)$ for $i<j$.

The model can be written as
\begin{equation}\label{bftl_nlane}
\begin{cases}
\dot x_n=v_n \\
\dot v_n = \alpha (V_j(\Delta x^j_n) - v_n)+\beta \dfrac{\Delta v^j_n}{(\Delta x^j_n)^2} & n\in I_j \\

\text{ + lane changing conditions}
\end{cases} \qquad \text{for } j=1,\dots,J.
\end{equation}

We adopt the lane changing conditions as in \eqref{rules}. Note that, except for the cases $j=1$ or $j=J$, if $j>2$ a vehicle might have the possibility to changes from lane $j$ to lane $j-1$ or  from lane $j$ to lane $j+1$. Consequently if both changes are possible we choose the most advantageous one in terms of acceleration.

As we done for the two-lane model we can define the steady state of model \eqref{bftl_nlane} in which all lane are at the equilibrium and lane changes do not occur. This is provided for the values of the headways \begin{equation}
\bar{h}_1,\dots,\bar{h}_J
\end{equation} that verify the condition
\begin{equation}\label{catenavel}
V_1(\bar{h}_1)=\cdots=V_J(\bar{h}_J).
\end{equation} In order to find this equilibrium we require also that the equilibrium velocity defined in \eqref{catenavel} must be smaller than the value $V_1^{max}$, that is the maximum velocity value allowed in the slower lane ($j=1$).

\subsection{An example with three lanes}
Let us consider a three-lane road ($J=3$). The steady state is given by the three values of the headways $\bar{h}_1,\bar{h}_2,\bar{h}_3$ such that  $V^{eq}:=V_1(\bar{h}_1)=V_2(\bar{h}_2)=V_3(\bar{h}_3)$. Using the same previous techniques can be show that with these conditions no lane changes occur and the system remains at the equilibrium.

We are now interested to add a perturbation in the middle lane and to study the possibility of lane changing. More specifically let us consider the initial conditions
\begin{equation}\label{dato3corsie}
\begin{array}{ll}
 \forall n\in I_1 
\begin{cases}
x_n(0)\,\,\text{equally spaced with distance}\,\,\bar{h}_1 \\
v_n(0)=V^{eq}
\end{cases}\\
\forall n\in I_2
\begin{cases}
x_n(0)\,\,\text{equally spaced with distance}\,\,\bar{h}_2+\eps \\
v_n(0)=V_2(\bar{h}_2+\eps)
\end{cases} \\
\forall n\in I_3
\begin{cases}
x_n(0)\,\,\text{equally spaced with distance}\,\,\bar{h}_3 \\
v_n(0)=V^{eq}.
\end{cases} 
\end{array}
\end{equation}

We can observe that the system is comparable to two subsystems: lane 1 - lane 2 and lane 2 - lane 3 where the lane changes are regulated by the thresholds in Table \ref{tabella1}. More specifically for the subsystem lane 2 - lane 3 we consider the case of a perturbation in the slow lane (first row of the table) while for the subsystem lane 1 - lane 2 we refer to the case of a perturbation in the fast lane (second row of the table). We add to this framework the possibility of choosing the best advantageous change for a vehicle in the middle lane that might have two possibilities for change lane.
The thresholds that enable lane changes can be obtained from Table \ref{tabella1} with the appropriate modifications. We have
\begin{table}[H]
\centering
\begin{tabular}{ |c|c|c|c|  }
\hline
 & $1\to 2$ \& $3\to 2$ & $2 \to 1$ & $2 \to 3$\\
\hline 
\rule{0mm}{1cm}
pert. $\eps$ in lane 2 & $\eps>d_s>0$ & $ \eps<\dfrac{V_1(\bar{h}_1-d_s)-V_1(\bar{h}_1)}{\left(1+\frac{\gamma}{(\bar{h}_1-d_s)^2}\right) V_2'(\bar{h}_2)}<0$ &  $\eps<\dfrac{V_3(\bar{h}_3-d_s)-V_3(\bar{h}_3)}{\left(1+\frac{\gamma}{(\bar{h}_3-d_s)^2}\right) V_2'(\bar{h}_2)}<0$\\
\hline
\end{tabular}
\caption{\label{tabella2}Thresholds and perturbations.}
\end{table}

\begin{figure}[H]
\centering
\includegraphics[scale=0.6]{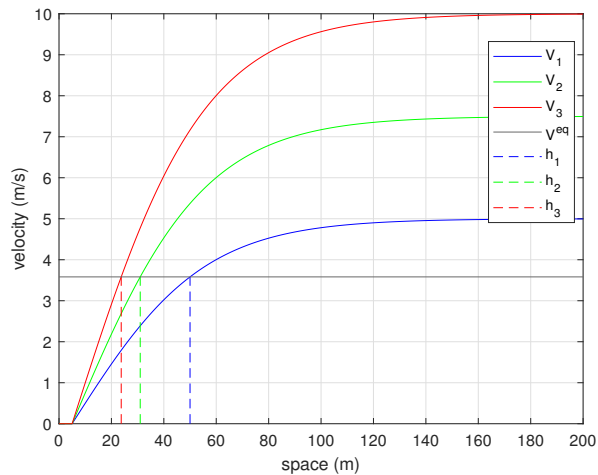}
\caption{Desired velocity functions.}\label{3vel}
\end{figure}

Here we propose a numerical example with a three-lane road, using the Runge Kutta 5 method. Consider the velocity function $V_1(h)$ as in \eqref{vel2corsie} and define $V_2(h)=\frac{3}{2}V_1(h)$ and $V_3(h)=2V_1(h)$. From the value $\bar{h}_1=50$ m we obtain that $\bar{h}_2=31$ m, $\bar{h}_3=23.7$ m and $V^{eq}=3.58$ m/s as shown in Fig. \ref{3vel}. The corresponding number of vehicles are: $\bar{N}_1=30$, $\bar{N}_2=48$, $\bar{N}_3=63$.
In order to add a perturbation in lane 2 we find the values of the perturbation that allow lane changing. From Table \ref{tabella2} we obtain: $\eps>5$ m for lane changes from lanes 1 and 3 to lane 2, $\eps<-2.25$ m for lane changes from lane 2 to lane 1 and $\eps<-7.74$ m for lane changes from lane 2 to lane 3.
In the following numerical tests we use the initial conditions \eqref{dato3corsie} with $\eps=-2.68$ m in the test (a) and with $\eps=-7.91$ m in the test (b). We can observe that in the test (a) the perturbation has produced lane changes from lane 2 to lane 1 while in the test (b) lane changes from lane 2 to lane 3 occurred.  
\begin{figure}[H]
\centering
\includegraphics[scale=0.35]{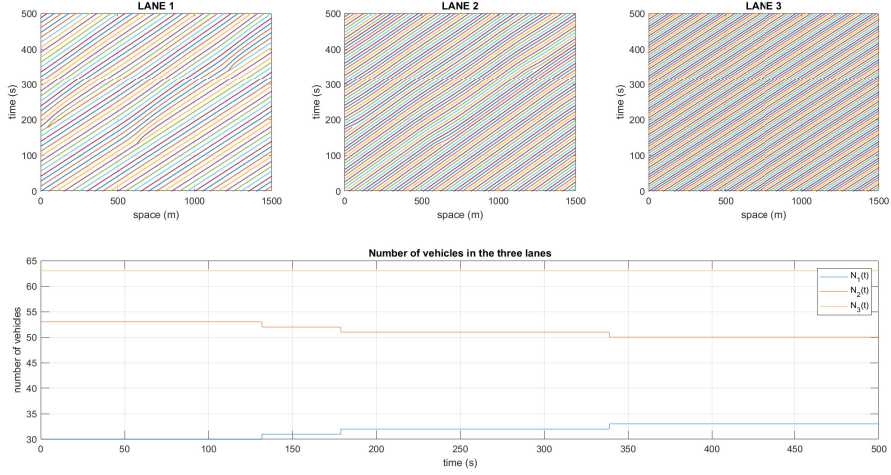}
\caption{Test (a) - Top: vehicle trajectories in the three lanes. Bottom: number of vehicles versus time.}\label{cambicorsia31}
\end{figure} 

\begin{figure}[H]
\centering
\includegraphics[scale=0.35]{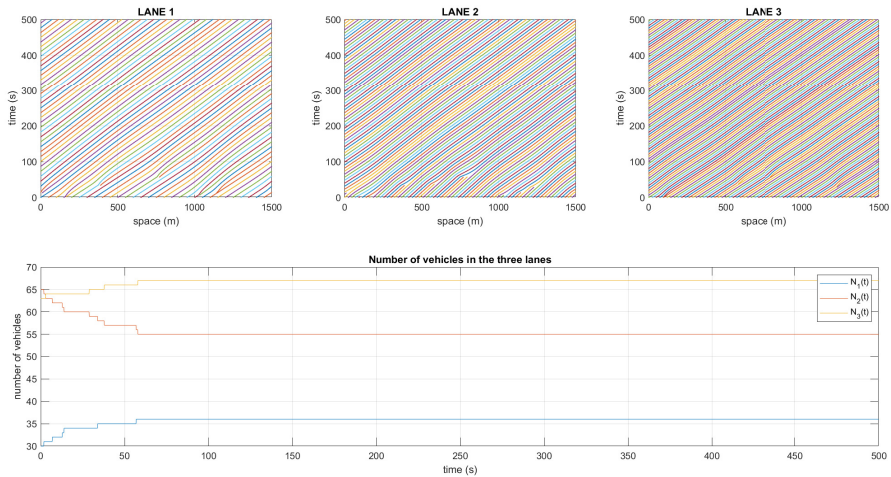}
\caption{Test (b) - Top: vehicle trajectories in the three lanes. Bottom: number of vehicles versus time.}\label{cambicorsia32}
\end{figure}

In the following test we show an example of instability, comparing the results with the test in section \ref{testsg}. Let us consider the function $V_1(h)$ as in the aforementioned test, and define $V_2(h)=\frac{3}{2}V_1(h)$ and $V_3(h)=2V_1(h)$. We consider $N_1(0)=N_3(0)=90$ and $N_2(0)=0$ with initial conditions with random perturbations $r_n$.
\begin{equation}
\begin{array}{ll}
 \forall n\in I_1 
\begin{cases}
x_n(0)-x_{n-1}(0)=\frac{L}{N_1(0)}+r_n \\
v_n(0)=V_1(\frac{L}{N_1(0)})
\end{cases}\\
\forall n\in I_3
\begin{cases}
x_n(0)-x_{n-1}(0)=\frac{L}{N_3(0)} \\
v_n(0)=V_2(\frac{L}{N_3(0)})
\end{cases} 
\end{array}
\end{equation}

\begin{figure}[H]
\centering
\includegraphics[scale=0.35]{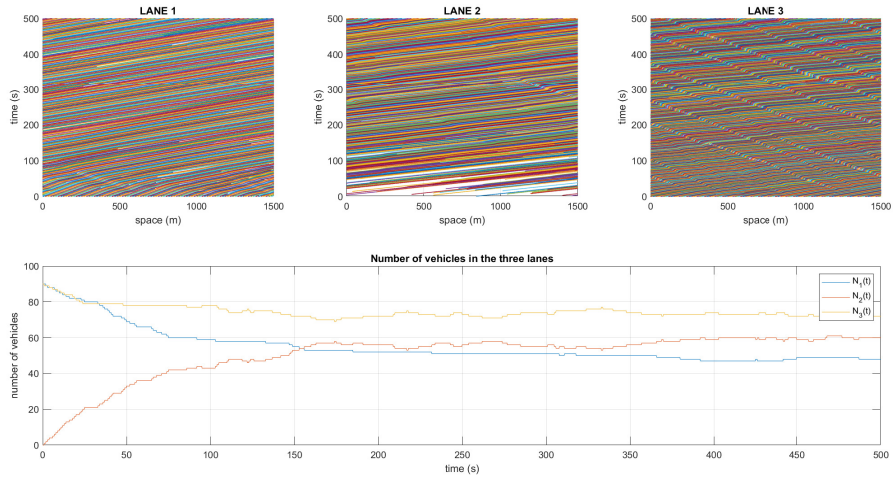}
\caption{ Top: vehicle trajectories in the three lanes. Bottom: number of vehicles versus time.}\label{cambicorsia3corsie}
\end{figure}

From Fig. \ref{cambicorsia3corsie} we can see that lane 1 gradually empties into lane 2. Due to frequent lane changes, more pronounced  {stop \& go} waves occur in fast lanes, while slow lane tends to stabilize. In test  \ref{testsg} we recall that the instabilities were evident in both lanes.

\section{Conclusions}  {
In this paper we have studied a microscopic model \eqref{bftl} for lane changing proposing simple lane changing rules. We have computed global steady states and we have investigated the linear stability of such solutions. The global steady state of the multi-lane model is parametrized by the total number $N$ of vehicles in the road. All lanes are coupled by the lane changing conditions, and the equilibrium is reached only when the crowding of each single lane is such that no lane changing is convenient anymore. At that point the system can reach the equilibrium lane by lane. We have proved that the model for the single-lane case has a larger stability region than the model \eqref{bando_classico}. In the multi-lane case we have proved that is possible to determine conditions on perturbations in which the equilibrium of the steady state is preserved and lane changing does not occur. We plan to derive a macroscopic version of this model where each lane would be described by its own equation and the lane changes  would appear as source terms for the macroscopic equations. This study can be useful in applications for instance in the design of velocity profiles to minimize lane changes  in order to avoid jams and car accidents. }

\subsection*{Acknowledgement}
Both authors are members of the INdAM Research group GNCS. This work was supported in part by Progetto di Ateneo 2019, n. 1622397 and 2020 n. 2023082 (Sapienza - Universit\`{a} di Roma), and PRIN 2017KKJP4X.

\nocite{*}
\printbibliography[heading=bibintoc]

\end{document}